\documentclass[11pt]{amsart}

\usepackage{amssymb,amsmath,amscd}
\usepackage{bbm}
\usepackage{a4wide}

\newtheorem{theorem}{Theorem}
\newtheorem{prop}{Proposition}
\newtheorem{lemma}{Lemma}
\newtheorem{coro}{Corollary}
\newtheorem{fact}{Fact}

\theoremstyle{definition}
\newtheorem{remark}{Remark}

\newcommand{\ts}{\hspace{0.5pt}}
\newcommand{\HH}{\mathbb{H}\ts}
\newcommand{\RR}{\mathbb{R}}
\newcommand{\CC}{\mathbb{C}}
\newcommand{\QQ}{\mathbb{Q}\,}
\newcommand{\ZZ}{\mathbb{Z}}
\newcommand{\NN}{\mathbb{N}\ts}

\newcommand{\II}{\mathbb I}
\newcommand{\JJ}{\mathbb J}
\newcommand{\KK}{\mathbb K\ts}

\newcommand{\LL}{\mathcal{L}}
\newcommand{\OO}{\mathcal{O}}
\newcommand{\oo}{{\scriptstyle \mathcal{O}}}

\newcommand{\gG}{\varGamma}
\newcommand{\one}{\mathbbm 1}
\newcommand{\bt}{{\scriptscriptstyle \bullet}}

\newcommand{\pl}{\! + \!}
\newcommand{\mi}{\! - \!}

\newcommand{\ii}{\mathrm{i}}
\newcommand{\jj}{\mathrm{j}}
\newcommand{\kk}{\mathrm{k}}

\newcommand{\R}{\operatorname{Re}}
\newcommand{\I}{\operatorname{Im}}
\newcommand{\N}{\operatorname{N}}
\newcommand{\tr}{\operatorname{tr}}
\newcommand{\nr}{\operatorname{nr}}
\newcommand{\cont}{\operatorname{cont}}
\newcommand{\lden}{\operatorname{lden}}

\begin{document}

\title[Coincidence site modules]
{Coincidence site modules in $3$-space}

\author{Michael Baake}
\address{Fakult\"at f\"ur Mathematik, Univ.\ Bielefeld,
Box 100131, 33501 Bielefeld, Germany}
\email{mbaake@math.uni-bielefeld.de}
\urladdr{http://www.math.uni-bielefeld.de/baake/}

\author{Peter Pleasants}
\address{Department of Mathematics, University of Queensland,
Brisbane, QLD 4072, Australia}
\email{pabp@maths.uq.edu.au}
\urladdr{http://www.maths.uq.edu.au/}

\author{Ulf Rehmann}
\address{Fakult\"at f\"ur Mathematik, Univ.\ Bielefeld,
Box 100131, 33501 Bielefeld, Germany}
\email{rehmann@math.uni-bielefeld.de}
\urladdr{http://www.math.uni-bielefeld.de/rehmann/}

\begin{abstract}
  The coincidence site lattice (CSL) problem and its generalization to
  $\ZZ$-modules in Euclidean $3$-space is revisited, and various
  results and conjectures are proved in a unified way, by using 
  maximal orders in quaternion algebras of class number 1
  over real algebraic number fields.
\end{abstract}

\maketitle

\vspace*{-5mm}\thispagestyle{empty}
\section{Introduction}

A \emph{lattice} in Euclidean space $\RR^d$ is a $\ZZ$-module of rank
$d$ whose $\RR$-span is $\RR^d$.  Two lattices $\gG$ and $\gG'$ in
$\RR^d$ are called \emph{commensurate}, denoted by $\gG\sim\gG'$, if
their intersection $\gG\cap\gG'$ has finite index both in $\gG$ and in
$\gG'$. This definition also applies to the commensurability of
$\ZZ$-modules in $\RR^d$ of rank $m>d$.  The \emph{commensurator} of a
lattice $\gG$ is the subgroup in $\mathrm{GL} (d,\RR)$ defined by
\[
   \mathrm{comm} (\gG) \; := \;
   \{ M\in\mathrm{GL} (d,\RR) \mid \gG\sim M\gG \}\, .
\]
It is not difficult to see that $\mathrm{comm} (\gG)\simeq
\mathrm{GL} (d,\QQ)$, via the map that takes the standard
basis of $\RR^d$ to a lattice basis.
More generally, for a subgroup $G\subset\mathrm{GL} (d,\RR)$, one
defines the commensurator of $\gG$ relative to $G$ in an analogous
way, denoted by $\mathrm{comm}^{}_{G} (\gG)$. Any commensurator of
this form is a function of the equivalence class defined by $\sim$,
i.e., commensurate lattices possess the same commensurators. One can
extend the definition to the affine group and subgroups thereof, but
that is outside our scope here.

Of particular interest in crystallography is the commensurator of a
lattice relative to the group of isometries, $\mathrm{O} (d,\RR)$, or
to its rotation subgroup, $\mathrm{SO} (d,\RR)$. This has been
investigated in great detail in connection with the classification of
grain boundaries and twins in materials science, compare
\cite{Boll,G}. The connection with the algebraic problem we are
interested in here emerges via the so-called coincidence site lattices
(CSLs), see \cite{BP,B,PBR} and references therein. These are finite
index sublattices of a given lattice $\gG$ that are intersections of
the form $\gG\cap R\gG$ with an isometry $R$ from the commensurator.  
In this context, the notation
\[
   \mathrm{(S)OC} (\gG) \; := \;
   \mathrm{comm}^{}_{\mathrm{(S)O} (d,\RR)} (\gG)
\]
has proved useful, together with the index function
$\varSigma\!:\,\mathrm{OC} (\gG) \longrightarrow \NN$ defined by
\[
    \varSigma (R) \; := \; [\gG : (\gG\cap R\gG ) ] \, .
\]
The integer $\varSigma (R)$ is the \emph{coincidence index}
of the isometry $R$, and the image 
$\varSigma \bigl(\mathrm{OC} (\gG)\bigr)$ is
called the elementary (or simple) \emph{coincidence spectrum} of
$\gG$. It refers to indices that emerge from single intersections
only.

More generally, one considers multiple intersections of the form
$\gG\cap R_1 \gG\cap\ldots\cap R_{\ell} \gG$, with each $R_i \in
\mathrm{OC} (\gG)$, and defines
$\varSigma (R_1,\ldots,R_{\ell})$ as the index of this intersection
in $\gG$. Clearly, the corresponding multiple spectra satisfy
\[
   \varSigma \bigl( \mathrm{OC} (\gG)^m \bigr) \; \subseteq \;
   \varSigma \bigl( \mathrm{OC} (\gG)^{n} \bigr)
   \quad\mbox{for}\quad m\le n ,
\]
and $\varSigma_{\rm tot} := \bigcup_{m\ge 1} \varSigma \bigl(
\mathrm{OC} (\gG)^m \bigr)$ is called the total (or complete)
coincidence spectrum of $\gG$.  In several important examples, compare
\cite{BG,Z3}, these spectra are equal or stabilize, in the sense that
the total spectrum is reached after finitely many lattice intersections.
Usually, there are more multiple CSLs than simple ones with a given
index, though this quantity also stabilizes in the examples mentioned.

Following the discovery of quasicrystals in the early 1980s, there was a
clear need for generalizing these concepts to accommodate aperiodic
situations.  A very natural approach consists in extending the setting
from lattices to finitely generated $\ZZ$-modules that are embedded in
$\RR^d$, such as rings of cyclotomic integers in the plane
\cite{PBR,BG} or similar sets in higher dimensions \cite{B,Z1,Z2}. In
general, such sets are dense, but commensurability and the
commensurator are still well defined, all on the basis of
group-subgroup indices, see \cite{B} for details.  This is the
approach we shall use here, too. For a given $\ZZ$-module $\gG$, it
leads to the determination of its elementary coincidence spectrum and
the corresponding coincidence site modules (CSMs), classified
according to their index in $\gG$.  In this paper, we focus on single
intersections, as they permit a rather general solution with algebraic
methods based on maximal orders in quaternion algebras over real
algebraic number fields. Single intersections are also the most
important ones for the applications in physics, as they refer to the
grain boundaries of twins, which are more abundant than multiple
junctions.

The paper is written in a self-contained way, in order to also reach
readers from the applied sciences. For this reason, we have not aimed
at maximal generality, and we have tried to give simple proofs, rather
than refer to general abstract results (though the latter are also
mentioned at several places). The full machinery will be developed
in a forthcoming publication \cite{BPnew}, together with further
examples that are beyond our scope here.

The material is organized as follows.  Section~\ref{example} reviews
the situation of the body centred cubic lattice as a motivating
example, hinting at the underlying connection to quaternion
arithmetic. The connection with quaternion algebra is summarized in
Section~\ref{gen-quat}, followed by the core of the paper in
Sections~\ref{common} -- \ref{general}. Here, we develop the precise
connection to reduced one-sided ideals in maximal orders and the
general solution of the coincidence problem in the case that both the
quaternion algebra and the base field have class number $1$.
Section~\ref{specific} applies the general findings to three specific
cases of interest, which are related to the Hurwitz, the icosian and
the octahedral rings, followed by some open questions.

\section{Example:\ The body centred cubic lattice in $\RR^3$}
\label{example}

Let $\ZZ^3$ be the standard primitive cubic lattice in $3$-space,
and consider $\gG_{\rm bcc} := \ZZ^3 \cup (v + \ZZ^3)$, where
$v=\frac{1}{2} (1,1,1)^t$. It is well known \cite{G,BP} that
\[
   \mathrm{(S)OC} (\gG_{\rm bcc}) \; = \;
   \mathrm{(S)O} (3,\QQ) \, ,
\]
which is clear from $\gG_{\rm bcc} \sim \ZZ^3$ and the corresponding
statement for $\ZZ^3$, compare \cite{B,BP}. From now on, we shall
restrict ourselves to rotations, i.e., to orientation preserving
isometries. In $3$-space, the orientation reversing isometries are of
the form $-R$ with $R$ a rotation. As $- \one^{}_{3}$ commutes with
all isometries and is a symmetry of $\gG_{\rm bcc}$ (and, in fact, of
all lattices and $\ZZ$-modules $\gG$ in $3$-space), it is clear that
the consideration of $\mathrm{SOC} (\gG_{\rm bcc})$ (or of
$\mathrm{SOC} (\gG)$) is sufficient to analyze the corresponding
coincidence problem completely. Note also that $\mathrm{SO} (3,\QQ)$
is a dense subgroup of $\mathrm{SO} (3,\RR)$.

The main interest lies in the classification of the CSLs
$\gG_{\rm bcc} \cap R \gG_{\rm bcc}$ that emerge from
$R\in\mathrm{SO} (3,\QQ)$, grouped according to their index
\[
   \varSigma(R) \; := \; 
   [\gG_{\rm bcc} : (\gG_{\rm bcc} \cap R \gG_{\rm bcc} )]\, .
\]
The answer is well known \cite{G,B} and can be formulated either
explicitly in terms of lattices and their bases \cite{G} or in terms
of quaternions \cite{B,Z1}.  In particular, the elementary coincidence
spectrum is the set of odd integers,
\begin{equation} \label{spec1}
   \varSigma \bigl(\mathrm{SOC} (\gG_{\rm bcc})\bigr) \; = \;
   2\NN_{0} + 1 \; = \; \{1,3,5,\ldots\}\, .
\end{equation}

An interesting quantity is the number of CSLs of $\gG_{\rm bcc}$
of index $m$, denoted by $f(m)$. It turns out to be a multiplicative
arithmetic function, i.e., $f(mn)=f(m)f(n)$ for $m,n$ coprime, and
is thus specified by its values at prime powers (other than $1$):
$f(2^r)=0$ and $f(p^r)=(p+1) p^{r-1}$, for odd primes $p$.  This
results in the Euler factors
\[
   E_p (s) \; = \; \sum_{r\ge 0} \frac{f(p^r)}{p^{rs}}
   \; = \; \frac{1+p^{-s}}{1-p^{1-s}}
   \; = \; \frac{1-p^{-2s}}{(1-p^{-s})(1-p^{1-s})}
\]
for odd primes, and in $E^{}_{2} (s)= 1$.
The corresponding Dirichlet series generating function reads
\begin{equation} \label{bcc-diri}
  \begin{split}
   \varPhi_{\rm cub} (s) & \; = \; \sum_{m=1}^{\infty}
   \frac{f(m)}{m^s} \; = \; \prod_{p\in\mathcal{P}} E_p (s)
   \; = \; \frac{1-2^{1-s}}{1+2^{-s}}\cdot
   \frac{\zeta(s)\,\zeta(s-1)}{\zeta(2s)}  \\
   & \; = \;
   1 + \frac{4}{3^s} + \frac{6}{5^s} + \frac{8}{7^s}
   + \frac{12}{9^s} + \frac{12}{11^s} + \frac{14}{13^s} 
   + \frac{24}{15^s} + \frac{18}{17^s} + \frac{20}{19^s} 
   + \ldots  \, ,
   \end{split}
\end{equation}
where $\mathcal{P} = \{2,3,5,7,11,\ldots \}$ is the set of rational
primes and $\zeta(s)=\sum_{n\ge 1} n^{-s}$ is Riemann's zeta function.
As we shall see later, this generating function applies to all three
cubic lattices (primitive, face centred and body centred), which is
why we have called it $\varPhi_{\rm cub}$.

This Dirichlet series permits the determination of asymptotic
properties of $f(m)$, e.g., by means of Delange's theorem (though
we do not need the full strength of it), see \cite[Appendix]{BM} for a
formulation adapted to the situation at hand. The result is formulated
via the corresponding summatory function and reads
\[
   F(x) \; := \; \sum_{m\le x} f(m) \; \sim \;
   \frac{3}{\pi^2}\, x^2  \; = \; 
   \frac{1}{\zeta(2)}\,\frac{x^2}{2}\, , \qquad \mbox{as $x\to\infty$},
\]
where $1/\zeta(2)=6/\pi^2$ is the residue of $\varPhi_{\rm cub} (s)$
at its right-most singularity, a simple pole at $s=2$. This shows
that $f(m)$ increases linearly on average and, following
\cite{HW}, see Chapter 18.2 and the footnote on p.~266 of it, one also
says that $6\ts m/\pi^2$ is the average size of the arithmetic function
$f(m)$.

Moreover, as observed in \cite{B}, one can express \eqref{bcc-diri} in
terms of the zeta function of the Hurwitz ring $\JJ$ of integer
quaternions \cite{H} in the quaternion algebra $\HH(\QQ)$ as
\begin{equation}  \label{J-diri}
   \varPhi_{\rm cub} (s) \; = \; \frac{1}{1+2^{-s}}\cdot
   \frac{\zeta^{}_{\JJ}(s/2)}{\zeta(2s)} \; = \;
   \frac{\zeta^{}_{\JJ}(s/2)}{\zeta^{}_{\JJ.\JJ}(s/2)}\, ,
\end{equation}
where $\zeta^{}_{\JJ}(s)$ and $\zeta^{}_{\JJ.\JJ}(s)$ are the Dirichlet
series generating functions for the one- and two-sided ideals of $\JJ$,
see Section~\ref{Hurwitz-results} for details. This indicates that one
should be able to re-derive the result using the arithmetic structure
of the ring $\JJ$ and to generalize it to other situations. That is
what we achieve in this article, with special focus on results and
conjectures made in \cite{BP,B}.  For first results on multiple CSLs
of $\gG_{\rm bcc}$, we refer to \cite{Z2,Z3}.

\section{Quaternion algebra and orthogonal groups}
\label{gen-quat}

Let us first recall some basic results about quaternions and introduce
some notation, see \cite[Sec.~57]{O} for background material. Let $K$
be a real algebraic number field of finite degree over $\QQ$. We
consider the corresponding quaternion algebra,
\[
   \HH(K) \; = \; K + \ii K + \jj K + \kk K ,
\]
which is a skew field, otherwise known as a division algebra.  A
convenient form of the defining relations for the generating elements
$1$ (implicit in the above representation) and $\ii,\jj,\kk$ is the
original one due to Hamilton \cite{KR},
\[
   \ii^2 \, = \, \jj^2 \, = \, \kk^2 \, = \, \ii\jj\kk 
   \, = \, -1,
\]
together with the requirement that $K$ is central in $\HH(K)$.  We note that
$K$ is the exact centre of the algebra $\HH(K)$.

A single quaternion $x = x_0 + \ii x_1 +\jj x_2 + \kk x_3$ is
sometimes also written as $x=(x_0,x_1,x_2,x_3)$, i.e., as a row
vector, while $x^t$ denotes its column counterpart. The following
result is standard, compare \cite{KR} and references therein.

\begin{fact}
$\HH(K)$ is canonically equipped with the following structural maps
and properties. 
\begin{itemize}
\item[(i)] the \emph{conjugation}$:$ $x \mapsto \bar{x} := x^{}_0 -\ii
  x^{}_1-\jj x^{}_2-\kk x^{}_3$, which is an anti-automorphism on the
  $K$-algebra\/ $\HH(K)$, i.e., one has $\overline{xy} = \bar{y}\,
  \bar{x}$ for all $x,y \in \HH(K) ;$
\item[(ii)]
  the \emph{reduced norm} and \emph{reduced trace}$:$
  $\nr,\, \tr \! :\, \HH(K) \longrightarrow K$, defined by 
\[
   \nr(x) \, = \,  x\bar{x} \, = \,  \mbox{$\sum_i$}\;  x_i^2
   \quad \mbox{and} \quad 
   \tr(x) \, = \,  x + \bar{x} \, = \, 2 x_0  .
\]
  They satisfy the identity\/  $\nr(x+y) = \nr(x) + \nr(y) + \tr(x\bar{y})$,
  which makes\/ $\nr$ a quadratic form on the vector space\/ $\HH(K) ;$
\item[(iii)] the Euclidean scalar product associated to the reduced
  norm, $\langle x, y \rangle := \frac{1}{2}\tr(x\bar{y}) = \sum_i x_i
  y_i $, which defines the Euclidean metric on\/ $\HH(K) \cong K^4$
  and makes $\{1,\ii,\jj,\kk\}$ an orthonormal basis of\/ $\HH(K)$
  $($a so-called $4$-bein\/$);$
\item[(iv)] 
 an orthogonal decomposition $\HH(K) = \HH_+ \oplus \HH_{-}$
 into eigenspaces of the conjugation map, where
\[
  \begin{aligned}
 \HH_+  & \, := \, \{x\in\HH(K)\;|\; \bar{x} = x\} =  K
         \quad \mbox{and} \\
 \HH_-  & \, := \, \{x\in\HH(K)\;|\; \bar{x} = -x\} 
 \, = \, \ii K \oplus \jj K \oplus \kk K .
   \end{aligned}
\]
 The elements of\/ $\HH_+$ are the \emph{scalar quaternions} and
 those of\/ $\HH_-$ are the \emph{pure quaternions}.
 For $x \in \HH(K)$, we call the corresponding components\/ 
 $\R(x):= \frac{1}{2} (x + \bar{x})$ the \emph{real part}
 and\/ $\I (x):= \frac{1}{2} (x - \bar{x}) $ the \emph{imaginary part} of $q$. 
\qed
\end{itemize}
\end{fact}

For $x \in \HH(K)$, the (reduced) characteristic polynomial (and the minimal
polynomial for $x \not\in K$) is just $X^2 - \tr(x)X + \nr(x)$. It is
preserved under each algebra automorphism of $\HH(K)$, and the same is true
for conjugation, orthogonality, the metric, and the decomposition into scalar
and pure part.  Therefore, any inner automorphism $i_q \! : \, x \mapsto q x
q^{-1}$ (for $x\in \HH(K)$ and $q \in \HH(K)^\bt := \HH(K)\setminus \{0\}$ )
leaves both subspaces $\HH_+$ and $\HH_-$ invariant and induces an isometry on
$\HH_- \cong K^3$.  Hence, $q \mapsto i_q$ defines a homomorphism $\HH(K)^\bt
\longrightarrow \mathrm{SO}(\HH_-)$ where $\mathrm{SO} (\HH_-) \simeq 
\mathrm{SO} (3,K)$.

The map $i_q$ fixes the elements of the space $K+qK$, which, for $q
\not\in K$, intersects $\HH_-$ in the line $\I(q)K$, whose
orthocomplement $W$ in $\HH_-$ is a plane.  Hence, $\HH_- = \I(q)K
\oplus W$.  For $w \in W$, we obtain $0 = 2 \langle q ,  w \rangle  =
\tr(q\bar{w}) = q\bar{w} + w\bar{q} = -qw + w\bar{q}$. This implies
the relation $qwq^{-1} = q\bar{q}^{-1}w$ which shows that $i_q$,
restricted to $\HH_-$, is a rotation with axis $\I(q)$ and 
rotation number $q\bar{q}^{\; -1}$ (whose real part is $\cos(\varphi)$).

If $q \in \HH_- \setminus \{0\}$, then $q\bar{q}^{\; -1} = -1$, hence the
eigenvalues of $-i_q$ on $\HH_-$ are $-1$ on the one-dimensional subspace
$qK$, and $+1$ on its two-dimensional complement $W$, i.e., $-i_q$ is the
reflection symmetry of $\HH_-$ with respect to $qK$.  Recall that any rotation
of the three-dimensional space $\HH_-$ is a product of two reflection
symmetries $-i_q,\; -i_{q'}$ for $q,q' \in \HH_-$, see \cite[Thm.~43.4,
p.~105]{O} for details in the generality needed. Consequently, the rotation is
induced by $i_{qq'}$, which proves the surjectivity of the homomorphism
$\HH(K)^\bt \longrightarrow {\rm SO}(\HH_-)$ defined by $q \mapsto i_q$.

The above argument also shows that, for any $q \in \HH(K)$, one has
$qx = xq$ if and only if $x \in K + qK$, i.e., $q$ or $x$ is central,
or $\I(q)$ and $\I(x)$ are collinear (or parallel).

The homomorphism $q\mapsto i_q$ can be made more explicit as follows:\ 
After identifying $K^3 = \ii K \oplus \jj K \oplus \kk
K = \HH_-$, the general setting permits rotation matrices in
$\mathrm{SO} (3,K)$ to be written via Cayley's parametrization (which
goes back to Euler for the case $K=\RR$, compare \cite{KR}) as
\begin{equation} \label{Cayley3}
  R(q) \, := \, R_q \, = \, \frac{1}{\nr(q)}
  \begin{pmatrix} \kappa^2 \pl \lambda^2 \mi \mu^2 \mi \nu^2 &
  -2\kappa\nu + 2\lambda\mu & 2\kappa\mu + 2\lambda\nu \\
  2\kappa\nu + 2\lambda\mu & \kappa^2 \mi \lambda^2 \pl
  \mu^2 \mi \nu^2 & -2\kappa\lambda \pl 2\mu\nu \\
  -2\kappa\mu + 2\lambda\nu & 2\kappa\lambda + 2\mu\nu &
  \kappa^2 \mi \lambda^2 \mi \mu^2 \pl \nu^2 \end{pmatrix},
\end{equation}
where $q=\kappa+\ii \lambda+\jj \mu+\kk \nu\in \HH(K)^\bt$, and where
$K^3$ is considered to be represented by the column vectors
$(x,y,z)^t$ with $x,y,z \in K$.

This explicit parametrization is useful in practice.  The rotation
axis (when $R_q$ is not the identity) is given by
$(\lambda,\mu,\nu)^t$, while the corresponding rotation angle
$\varphi$, using that the trace of $R(q)$ is $1+2\,\cos(\varphi)$ by
Euler's theorem, is determined by
\[
   \cos (\varphi) \; = \; \mathrm{Re} (q/\bar{q}) \; = \;
   \frac{\mathrm{Re} (q^2)}{\nr(q)} \; = \;
   \frac{\kappa^2 -\lambda^2 -\mu^2 -\nu^2}
   {\kappa^2 +\lambda^2 +\mu^2 +\nu^2} \, .
\]

\begin{fact} \label{Cayleymap}
  If $K$ is a real algebraic number field,
  every rotation matrix $M\in\mathrm{SO} (3,K)$ is
  of the form $M = R_q$
  for some $q\in\HH (K)^\bt$, i.e., the mapping\/
  $\HH(K)^\bt \longrightarrow\mathrm{SO} (3,K)$ 
  defined by the Cayley parametrization\/
  $q\mapsto R_q$ from\/ $\eqref{Cayley3}$ is onto.
\end{fact}

\begin{proof}
  This follows from the above arguments. Alternatively, one can use the
  parametrization \eqref{Cayley3} for $\mathrm{SO} (3) = \mathrm{SO} (3,\RR)$,
  where it is well-known since Euler and emerges from the canonical $2$ to $1$
  map from $\mathrm{SU} (2) \longrightarrow \mathrm{SO} (3)$, combined with
  the observation that $R(q)= R(\alpha q)$ for arbitrary non-zero
  $\alpha\in\RR$. One can then show that a parametrization of the subgroup
  $\mathrm{SO} (3,K)$ can be achieved with the quaternions restricted to
  $\HH(K)^\bt$.
\end{proof}

Looking at this matrix setting in $4$-space, noting that $K^4 \cong \HH(K)
= K \oplus \ii K \oplus \jj K \oplus \kk K$ and recalling our convention
to write quaternions as row vectors, one has the relation
\[
    \bigl(q{a}q^{-1}\bigr)^t \; = \; \begin{pmatrix}
    1 & 0 \\ 0 & R_q \end{pmatrix} {a}^t
\]
for arbitrary ${a}\in\HH(K)$ and $q\in\HH (K)^\bt$, which reflects the
orthogonal decomposition into real and imaginary parts mentioned
above. In particular, one also has
\begin{equation} \label{matrix1}
   \I (q{a}q^{-1}) \; = \; R_q \I ({a}),
\end{equation}
which will be crucial for our further development.

\section{Quaternion arithmetic}

Consider the real algebraic number field $K$ and let $\oo$ be the ring of
integers in $K$, compare \cite{BS} for general background material on
algebraic number fields.

\medskip
\noindent{\bf Assumption.} \emph{We shall assume throughout
this article that $K$ has class number\/ $1$.}
\medskip

The important consequence of this, for our purposes, is that $\oo$ is
a principal ideal domain (PID) and hence that every
\mbox{$\oo$-module} we encounter will be free and thus has a basis.
Also, our results permit a simpler formulation this way, though
many can be generalized to avoid this assumption.

Next, we need to know what `integral' and `integer' mean in $\HH(K)$.
A quaternion $q$ is called \emph{integral} (over $\oo$), when both
$\tr(q)$ and $\nr(q) =: |q|^2$ are in $\oo$, i.e., when it is a root
of a monic polynomial with coefficients in $\oo$.

The notation $\lvert q\rvert^2$ instead of $\nr (q)$ is quite common
and highlights the fact that the reduced norm coincides with the
squared Euclidean norm of the quaternion under the standard
identification of the generating elements with the Euclidean $4$-bein.
However, it has to be stressed here that, in spite of the notation,
$\lvert q\rvert^2$ is in general not a square in $K$.

Clearly, if $q$ is integral, then so is $\bar q$.  Note
that, because of the invariance of conjugation, this definition is
invariant under all automorphisms of $\HH(K)$.  However, due to
non-commutativity, the set of integral quaternions fails to be a ring.
This is the origin of various differences from the arithmetic of
commutative number fields.

Following \cite{V}, we call a finitely generated $\oo$-module whose
$K$-span is the whole of $\HH(K)$ an \emph{ideal}.  An \emph{order} is
an ideal that is also a ring containing $1$, see \cite[Sec.~8]{R} and
\cite[Chap.~1, Sec.~4]{V} for details.  It is called \emph{maximal}
when it is not contained in any larger order of $\HH (K)$. Note that
some authors restrict the term `ideal' to $\oo$-modules whose left and
right orders are maximal, because their properties are then more
similar to those of ideals in algebraic number fields.  Another term
in use is that of a (full) $\oo$-lattice, compare \cite{R,V}.
   
The following basic properties of an order $L$ are
important throughout the article.

\begin{fact} \label{basics}
All elements of an arbitrary order $L\subset \HH (K)$ are integral.
Moreover, one has $\bar{L}=L$ and $L\cap K = \oo$.
\end{fact}

\begin{proof}
  The first assertion follows from \cite[Lemme~I.4.1 and
  Prop.~I.4.2]{V} or \cite[Thms.~8.6, 9.3 and 1.14]{R}.  If $x\in L$,
  one has $\bar{x}=\mathrm{tr} (x) - x$, so $\bar{x}\in L$ since
  $\mathrm{tr} (x) \in\oo$, and $\bar{L}=L$. 
  
  Since $1\in L$ and $L$ is an $\oo$-module, $\oo\subseteq L\cap K$ is clear.
  Conversely, if $x \in L\cap K$, one has $\mathrm{tr} (x) = 2x \in \oo$ and
  $\mathrm{nr} (x) = x^2 \in\oo$, which implies that also $x\in\oo$.
\end{proof}

An element $a$ of an order $L$ is an \mbox{\emph{$L$-unit}} when
$a^{- 1}\in L$, and it is a consequence of Fact~\ref{basics} that
$a\in L$ is an \mbox{$L$-unit} if and only if $\nr(a)$ is a unit
of $\oo$. When it is clear which order is meant we shall omit the
prefix $L$ and simply speak of a \emph{unit}.  A \emph{left ideal}
of $L$ is an ideal $I$ with $aI\subseteq I$ for every $a\in L$ and
it is \emph{principal} if it has the form $Lx$ for some $x\in\HH(K)$.
Right ideals and principal right ideals of $L$ are defined similarly.
A \emph{two-sided} ideal of $L$ is an ideal that is simultaneously
a left ideal and a right ideal of $L$. 

Clearly,
\begin{equation} \label{lattice}
   \LL \; := \; \oo + \ii\oo + \jj\oo + \kk\oo
\end{equation}
is an order, as is $q \LL q^{-1}$ for any $q\in\HH(K)^\bt $; but it is
not maximal because $\frac{1}{2}(1+\ii+\jj+\kk)$ is integral and
appending it to $\LL$ generates a larger order.  For the applications,
we are mainly interested in maximal orders that contain $\LL$. In the
literature, such maximal orders are sometimes simply referred to as
rings of integers.

\begin{fact} \label{LversusO}
  If $q$ is an element of a maximal order $\OO$ that contains $\LL$,
  then $2q\in\LL$.
\end{fact}
\begin{proof}
  The $\oo$-basis of $\LL$ is also a $K$-basis of $\HH(K)$, so $q=a +
  \ii b + \jj c + \kk d$ with $a,b,c,d\in K$. Since $q$, $\ii q$, $\jj
  q$ and $\kk q$ must all be integral (as $\ii$, $\jj$ and $\kk$ are
  units of $\OO$), their reduced traces must be in $\oo$.  They are
  $2a$, $-2b$, $-2c$ and $-2d$, respectively, so that $2q\in\LL$.
\end{proof}

Maximal orders $\OO$ are certainly not unique, as the application of
inner automorphisms results in different maximal orders.  They may not
even be unique up to inner automorphism.  Orders related by an inner
automorphism are said to be of the same \emph{type}.  The \emph{class
  number} of a maximal order $\OO$ is the number of equivalence
classes of left ideals of $\OO$ under the equivalence relation of
right multiplication by elements of $\HH(K)$.  Since the conjugation
anti-automorphism takes left ideals to right ideals, right
multiplication to left multiplication and $\OO$ to itself,
interchanging the roles of `left' and `right' in this definition
results in the same value of the class number.

\begin{fact} \label{types-and-classes}
All maximal orders of\/ $\HH(K)$ have the same class number, 
in the sense just defined, which is known as the \emph{class number 
of\/ $\HH(K)$}. Moreover, if the class number is\/ $1$, all
maximal orders are mutual images of one another under inner
automorphims.
\end{fact}

\begin{proof}
This follows from Lemme~I.4.10 and Exercise~I.4.1 of \cite{V}. 
\end{proof}

When $\HH(K)$ has class number $1$, every left ideal of a maximal order
$\OO$ has the form $\OO q$ for some $q\in\OO^\bt:=\OO\setminus\{0\}$
(and every right ideal has the form $q\OO$).  In the case of $\HH(\QQ)$,
there is only one type of maximal order, and a unique maximal
order $\OO$ that contains $\LL$, see \cite{H}:
\begin{equation} \label{J-def}
    \OO \; = \; \JJ \; = \; 
    \big\langle 1,\ii,\jj, \tfrac{1}{2} (1+\ii+\jj+\kk)
    \big\rangle^{}_{\ZZ} \, .
\end{equation}
This is the so-called Hurwitz ring of integers \cite{H}. It is an
index $2$ extension of $\LL \simeq \ZZ^4$. In fact, seen as a lattice
in 4-space, it is the weight lattice $D^*_4$, the dual of the root
lattice $D_4$, the latter in standard representation with roots of
squared length $2$, compare \cite{CSl}.

Later on, we shall need some information about the relation between
$\OO$ and $q\OO q^{-1}$, for $q\in\OO$. It is immediate that
\begin{equation}  \label{cond-for-twosided}
   \OO \, = \, q\OO q^{-1} \quad\Longleftrightarrow\quad
   q\OO \, = \, \OO q \quad\Longleftrightarrow\quad
   \bar{q}\OO \, = \, \OO\bar{q}\, ,
\end{equation}
the last characterization following from $\overline{\OO}=\OO$.
When this occurs, $q\OO = \OO q$ is said to be a \emph{two-sided
  ideal} of $\OO$, and $q$ is said to generate a two-sided ideal. All
elements of $K$ commute with $\OO$, and thus generate two-sided
ideals. Another example we shall meet below is $1+\ii$, which
generates a two-sided ideal in $\JJ$, but not in our other explicit
examples of Section~\ref{specific}.

{}Finally, we observe the following counterpart of the relation
$L\cap K = \oo$ from Fact~\ref{basics}, where we use the set-valued
extension of $\R$, i.e.,  $\R (S) := \{ \R (x)\mid x\in S\}$.
\begin{fact} \label{real-part-of-O}
   If\/ $\HH (K)$ has class number $1$ and if
   $\OO$ is a maximal order in\/ $\HH (K)$,  $\OO$  has
   real part\/ $\R (\OO) = \tfrac{1}{2} \oo$. 
\end{fact}

\begin{proof}
  Assume first that $\OO$ is a maximal order that contains $\LL$ of
  \eqref{lattice} as well as the integral quaternion
  $\tfrac{1}{2}(1+\ii+\jj+\kk)$. It is clear that such a maximal order
  exists, and that $\tfrac{1}{2} \oo \subseteq \R (\OO)$ in this situation.
  On the other hand, $2q\in\LL$ for all $q\in\OO$ by
  Fact~\ref{LversusO} in this case, and we must also have $\R (\OO)
  \subseteq \tfrac{1}{2} \oo$, hence $\R (\OO) = \tfrac{1}{2} \oo$.
  By Corollaire~I.4.11 and Exercise~I.4.1 of \cite{V},
   all other maximal orders of $\HH(K)$ are images of this
   one under inner automorphisms, and, since the real part
   is invariant under inner automorphisms, the result follows.
\end{proof}

Below, we shall often meet intersections of two (not necessarily
distinct) maximal orders. These are also known as \emph{Eichler orders},
see \cite{V,St} for details.

\label{common}

In this and the following two sections, we develop some general
arguments without using specific properties of the maximal orders we
have in mind for the applications. To this end, let $K$ as above be a
real algebraic number field with class number $1$ whose ring of
integers $\oo$ is therefore a PID.  The class number property implies
that every $\oo$-module in $\HH (K)$ is free (and thus has a basis)
\cite[Ch.~3, Thm.~6.6]{AW}, since such modules are torsion-free.
In what follows, $L$ will be an
arbitrary order in $\HH(K)$ and $\OO$ a maximal order. Later on, we
shall need the assumption that $\OO$ has class number $1$, i.e., that
all ideals of $\OO$ are principal.  This condition will be formulated as
a property of $\HH (K)$ and always be mentioned explicitly.

\smallskip
Let us return to the parametrization of the rotations.
It is clear that the choice of $q$ in \eqref{Cayley3} is not unique.
In fact, we can always arrange $q\in\OO^\bt$ (by multiplying $q$ by a
suitable element of $\oo$).  An element $q\in\OO$ is called
\emph{$\OO$-primitive} if it is not divisible by any
non-unit of $\oo$. For an arbitrary $q\in \OO^\bt$, one can define
\begin{equation} \label{def-cont}
    \mathrm{cont}^{}_{\OO} (q) \; := \;
    \mathrm{lcm} \{ \alpha \in\oo\mid q/\alpha \in \OO \}\, ,
\end{equation}
which is called the \emph{$\OO$-content} of $q$, see \cite{BM} for
details on this concept in a more general setting.  It is well defined
up to units of $\oo$ due to our assumption that $\oo$ is a PID.  So,
$q$ is $\OO$-primitive when $\mathrm{cont}^{}_{\OO} (q)$ is a unit in
$\oo$. Alternatively, one can consider the content to be the principal
ideal generated by $\mathrm{cont}^{}_{\OO} (q)$, which takes care of
the units automatically.

More generally, when $\oo$ is a PID as here, primitivity relative to
any order $L$ of $\HH(K)$ can be defined in an analogous way: $q\in L$
is \mbox{\emph{$L$-primitive}} if it cannot be factorized as $q=\alpha
r$, with $r\in L$ and $\alpha$ a non-unit in $\oo$.  As with units, we
shall omit the prefix $L$ when it is clear which order is meant.  Also
the \emph{$L$-content}, $\cont^{}_{L} (q)$, for any $q\in L^\bt$ can be
defined analogously to the $\OO$-content, and is again unique up to
units of $\oo$. In particular, primitivity relative to the order $\LL$
of \eqref{lattice} is another useful concept. In \cite{H}, an element
${q}\in\LL$ that is $\LL$-primitive in this sense is often simply
called \emph{primitive}, and this property is equivalent to the
condition that $\gcd({q}):=\gcd(q^{}_0,q^{}_1,q^{}_2,q^{}_3)=1$, where
the $\gcd$ is defined as usual, and determined up to a unit of $\oo$.
More generally, one has $\cont^{}_{\LL} (q)=\gcd (q)$.

\begin{prop} \label{Cayley-para}
  Let $L$ be an arbitrary order of\/ $\HH(K)$ and let\/ $R\! : \,
  q\mapsto R(q)$ be Cayley's parametrization of\/
  $\eqref{Cayley3}$. Then, one has $R(L^\bt) = \mathrm{SO}
  (3,K)$. Moreover, any rotation matrix $M\in\mathrm{SO} (3,K)$ can be
  parametrized as\/ $M=R(q)$ with an \mbox{$L$-primitive} element\/
  $q\in L$.  This gives a multiple cover of\/ $\mathrm{SO} (3,K)$,
  where the remaining multiplicity is due to\/ $R(\varepsilon q)=R(q)$
  with\/ $\varepsilon$ a unit in\/ $\oo$.
\end{prop}

\begin{proof}
Since $R(\HH(K)^\bt)=\mathrm{SO} (3,K)$ by Fact~\ref{Cayleymap}, it is
clear that, given a matrix $M\in\mathrm{SO} (3,K)$, $q$ can be
chosen to satisfy $R(q)=M$ and to be an \mbox{$L$-primitive} quaternion 
at the same time, the latter condition due to the freedom to
multiply $q$ by a suitable non-zero element of $K$, which commutes
with $q$ and does not change the matrix $M$.

Clearly, $R(\varepsilon q)=R(q)$ for any unit $\varepsilon\in\oo$, and
$\varepsilon q$ is \mbox{$L$-primitive} when $q$ is.  In view of
Eq.~\eqref{matrix1}, and the fact that $\R(qaq^{-1})=\R(a)$, two
quaternions ${p},q$ parametrize the same matrix $M$ if and only if $
{p}^{-1}q$ is a central element of $\HH(K)$, so $q = \alpha{p}$ with
$\alpha\in K$. When ${p}$ and $q$ are both $L$-primitive, such a relation
can only hold when $\alpha$ is a unit in $\oo$.
\end{proof}

\begin{remark} \label{Cayley-rem}
     The relevance of Proposition~\ref{Cayley-para} originates in the
     observation that
\begin{equation}
     \mathrm{SOC} \bigl(\I (\LL) \bigr) \; = \;
     \mathrm{comm}^{}_{\mathrm{SO} (3,\RR)} 
         \bigl(\I (\LL) \bigr) \; = \;
     \mathrm{SO} (3,K) \, ,
\end{equation}
since $\I (\LL)$ is the $\oo$-span of the standard Euclidean $3$-bein
in $\RR^3$. Moreover, $\mathrm{SO} (3,K)$ is the $\mathrm{SOC}$-group
of $\I (L)$ for \emph{any} order $L$ of $\HH(K)$, because $\LL$ and $L$ (and
hence also $\I (\LL)$ and $\I (L)$) are commensurate.
\end{remark}

\section{Index calculations} \label{indices-section}

One of our main goals is to compute indices of CSMs. As we shall see
below, this will require working out indices of the form $[\I (L) : \I
(qL)]$ for elements $q\in L^\bt$. To achieve this in a simple and
systematic fashion, it is advantageous along the way, and essential
for our later proofs in this concrete approach, to work with indices
which take values in $\oo$ rather than in $\ZZ$, and we define these
now. In fact, we define these more general indices to be principal
ideals in $\oo$.  In this setting, the coset-counting indices of
interest can later be derived as the absolute norms of the
corresponding $\oo$-indices, see Fact~\ref{dets} below.

For a vector space $V$ over $K$ and a $K$-linear map $\phi$ of $V$
into itself, we use $\det^{}_K (\phi)$ to denote the usual determinant
of $\phi$. It is an element of $K$ and can be calculated using a
$K$-basis of $V$. Similarly, $\det^{}_{\QQ} (\phi)$ denotes the
determinant of $\phi$ when $V$ is regarded as a vector space over
$\QQ$ and $\phi$ as a $\QQ$-linear map. It is a rational number and
can be calculated using a $\QQ$-basis of $V$. In both cases, the
result is independent of the particular basis chosen.

If $J \subseteq I$ are full $\oo$-modules in $V$, we define the
\emph{$K$-index} of $J$ in $I$ as
\begin{equation} \label{det}
    [ I : J ]^{}_{K} \; := \; \mathrm{det}^{}_K (\phi)\,\oo  ,
\end{equation}
where $\phi$ is a linear map that takes an $\oo$-basis of $I$ to an
$\oo$-basis of $J$ and $\det^{}_K(\phi)\, \oo$ is the principal ideal
of $\oo$ generated by $\det^{}_K(\phi)$. The determinant
$\det^{}_K(\phi)$ is in $\oo$, since $J\subseteq I$, and changing the
bases of $I$ and $J$ multiplies it by a unit of $\oo$, so the ideal
$[I:J]^{}_K$ is independent of the chosen bases.  By considering the
inverse of $\phi$, it is clear that $I=J$ if and only if
$[I:J]^{}_K=\oo$.  The $K$-index is multiplicative: if $I_1 \supseteq
I_2 \supseteq I_3$ are full $\oo$-modules in $V$, one has
\begin{equation} \label{multi}
   [ I_1 : I_3 ]^{}_K \; = \;
   [ I_1 : I_2 ]^{}_K \, [ I_2 : I_3 ]^{}_K .
\end{equation}  
We also write $[I:J]=\lvert\det^{}_{\QQ} (\phi) \rvert$, a positive integer,
which is the index of $J$ in $I$ in the usual sense.

\begin{lemma}\label{[L:qL]}
   For any order $L$ in\/ $\HH(K)$ and any $q\in L^\bt$, one has
   $\; [L:qL]^{}_K = [L:Lq]^{}_K = \lvert q\rvert^4 \oo$.
\end{lemma}

\begin{proof}
  It is sufficient to prove the claim for right ideals $qL$.  The
  linear map $\phi_q \! : \, x\mapsto qx$ (the left regular
  representation of $q$ on $\HH(K)$) takes an \mbox{$\oo$-basis} of
  $L$ to an \mbox{$\oo$-basis} of $qL$, and a straightforward
  calculation using the basis $\{1,\ii,\jj,\kk\}$ of $\HH(K)$ gives
  $\det^{}_K (\phi_q) =\lvert q\rvert^4$.
\end{proof}

Our next result is an index formula that provides the link to our
problem in $3$-space by giving the indices of images under the mapping
$x\mapsto\I(x)$ and paves the way to \mbox{$3$-dimensional} analogues
of Lemma~\ref{[L:qL]}.

\begin{prop}\label{indices}
If $J\subseteq I$ are ideals in\/ $\HH(K)$, one has the index relations
\begin{align*}
    [I:J]^{}_K \; &= \;
    [I\cap \HH_+ : J\cap \HH_+]^{}_K \, [\I(I) : \I(J)]^{}_K 
    \quad \mbox{and} \\
    [I:J]^{}_K \; &= \;
    [\R(I) : \R(J)]^{}_K \, [I\cap \HH_- : J\cap \HH_-]^{}_K \,,
\end{align*}
where all five indices are principal ideals of $\oo$.
\end{prop}

\begin{proof}
  These identities are both instances of the fact that the density of
  a lattice (being the reciprocal of the absolute value of its
  determinant) is the density of an orthogonal projection of it
  multiplied by the density of the kernel of the projection.

For the first identity, choose generators $\alpha^{}_0$ and
$\beta^{}_0$ of $I\cap \HH_+$ and $J\cap \HH_+$ and extend
them to \mbox{$\oo$-bases}
\[
    \{\alpha^{}_0,\alpha^{}_1+u^{}_1,
      \alpha^{}_2+u^{}_2,\alpha^{}_3+u^{}_3\}
    \quad \mbox{and} \quad
    \{\beta^{}_0,\beta^{}_1+u^{}_1,
      \beta^{}_2+u^{}_2,\beta^{}_3+u^{}_3\}
\]
of $I$ and $J$; these operations being possible by
\cite[Thm.~6.16]{AW}, since $\alpha^{}_0$ and $\beta^{}_0$ are
primitive elements of $I$ and $J$ in the sense of \cite{AW}.
By construction, $\{\I (u_1^{}), \I (u_2^{}), \I (u_3^{})\}$ 
and $\{\I (v_1^{}), \I (v_2^{}), \I (v_3^{}) \}$
are then \mbox{$\oo$-bases} of $\I(I)$ and $\I(J)$. Note that,
due to our convention, each $\I (u_i)$ is the
transpose of the row vector $u_i$ without its first
coordinate, and similarly for the row vectors $v_j$.
Since $I\cap \HH_+$ and $J\cap \HH_+$ are \mbox{$1$-dimensional},
$[I\cap \HH_+:J\cap \HH_+]^{}_K$ is trivially evaluated as
$(\beta^{}_0/\alpha^{}_0)\oo$.  If
\[
   x \; = \; \mu^{}_0 \alpha^{}_0+
     \mu^{}_1 (\alpha^{}_1+u^{}_1)+
     \mu^{}_2 (\alpha^{}_2+u^{}_2)+
     \mu^{}_3 (\alpha^{}_3+u^{}_3)\, ,
\]
one has
\[
   \I(x) \; = \; \I \bigl(\mu^{}_1 u^{}_1 +
         \mu^{}_2 u^{}_2 + \mu^{}_3 u^{}_3 \bigr) ,
\]
so the $4\!\times\! 4$ matrix of the linear map taking the basis
of $I$ to the basis of $J$ has the form
\[
   M \; = \; \begin{pmatrix} \beta^{}_0 /\alpha^{}_0 & m \\
    0 & M_3 \end{pmatrix}^t ,
\]
where $M_3$ is the $3\!\times\! 3$ matrix of the linear
map taking the basis of $\I(I)$ to the basis of $\I(J)$.
Taking determinants of both sides gives
\[
   [I:J]^{}_K \, = \, \det (M)\,\oo \, = \,
   (\beta^{}_0 /\alpha^{}_0) \det (M_3)\,\oo \, = \,
   [I\cap \HH_+:J\cap \HH_+]^{}_K \, [\I(I):\I(J)]^{}_K.
\]

For the second identity, let $\alpha$ and $\beta$ be generators of the
$\oo$-ideals $\R(I)$ and $\R(J)$, choose $u^{}_0$ and $v^{}_0$ in
$\HH_-$ with $\alpha+u^{}_0\in I$ and $\beta+v^{}_0\in J$, and choose
bases $\{u_1,u_2,u_3\}$ and $\{v_1,v_2,v_3\}$ of $I\cap\HH_-$ and $J
\cap\HH_-$. Then, $\{\alpha+u_0,u_1,u_2,u_3\}$ and
$\{\beta+v_0,v_1,v_2,v_3\}$ are bases of $I$ and $J$, and the
$4\!\times\! 4$ matrix of the linear map taking the former to the
latter is
\[
   M \; = \; \begin{pmatrix} \beta/\alpha & m \\
    0 & M_3^t \end{pmatrix},
\]
where $M_3$ is the $3\!\times\! 3$ matrix of the linear mapping that
takes $\{\I (u_1^{}), \I (u_2^{}), \I (u_3^{})\}$ to $\{\I (v_1^{}),
\I (v_2^{}), \I (v_3^{}) \}$.  Hence
\[
   [I:J]^{}_K \, = \, \det (M)\,\oo \, = \,
   (\beta/\alpha) \det (M_3)\,\oo \, = \,
   [\R(I):\R(J)]^{}_K \, [I\cap \HH_-:J\cap \HH_-]^{}_K \, ,
\]
which completes the argument.
\end{proof}

Since $\HH_{+} \simeq K$ in our setting, we shall usually
write $I\cap K$ instead of $I\cap \HH_{+}$ from now on.

\begin{fact}\label{im-index}
   If $L' \subset L$ are orders in\/ $\HH(K)$, one has
   $\; [L:L']^{}_K=[\I(L):\I(L')]^{}_K$.
\end{fact}

\begin{proof}This follows from Proposition~\ref{indices},
since $L\cap K=L'\cap K=\oo$, by Fact~\ref{basics}.
\end{proof}

\begin{lemma}\label{qLcapK}
   If $L$ is an order of\/ $\HH(K)$ and $q\in L$ is
   $L$-primitive, one has
\[
    qL\cap K  \, = \, L\bar{q}\cap K \, = \,
    \lvert q\rvert^2\oo .
\]
\end{lemma}

\begin{proof}
  It is sufficient to prove $qL\cap K = \lvert q\rvert^2 \oo$, as the
  other claim follows by conjugation. Observe that the inclusions
  $\lvert q\rvert^2 \oo\subseteq qL\cap K$ and $qL\cap K\subseteq \oo$
  are clear.  It remains to show that $\alpha\in qL\cap K$ implies
  $\lvert q\rvert^2 \mid\alpha$. We have $\bar{q}\alpha =
  \alpha\bar{q}=\lvert q\rvert^2r$, where also $\bar{q}$ is
  $L$-primitive, i.e., $\mathrm{cont}^{}_{L} (\bar{q})$ is a unit in
  $\oo$. Consequently, we have $\lvert q\rvert^2 \mid 
  \mathrm{cont}^{}_{L} (\alpha \bar{q}) = \alpha$.
\end{proof}

\begin{lemma}\label{[ImL:ImqL]}
   If $L$ is an order of\/ $\HH(K)$ and $q\in L$ is $L$-primitive,
   one has the index formula
\[
   [\I(L):\I(qL)]^{}_K \, = \,
   [\I(L):\I(L\bar{q})]^{}_K
    \, = \, \lvert q\rvert^2 \oo .
\]
\end{lemma}
\begin{proof}
This follows from Proposition~\ref{indices} and Lemmas~\ref{[L:qL]}
and \ref{qLcapK}, since $[\oo:\lvert q\rvert^2\oo]^{}_K=\lvert q\rvert^2 
\oo$.
\end{proof}

In the three lemmas derived so far, and in various places below,
there are obvious additional identities, due to the freedom to
replace $q$ by $\bar{q}$.

To make the transition from $K$-indices to the usual indices of
modules, we need to introduce norms from $K$ to $\QQ$. Let $\N:=\lvert
\N_{K/\QQ}\rvert$ be the absolute norm from $K$ to $\QQ$, which is
given by $\N(\alpha)=\big\lvert\prod\sigma(\alpha)\big\rvert\in\QQ$,
where $\sigma$ runs through the homomorphisms of $K$ into $\CC$.  It
is a non-negative rational number for all $\alpha\in K$, and a
non-negative integer for $\alpha\in\oo$.

The absolute norm can also be defined for principal ideals of
$\oo$ by $\N(\alpha\oo)=\N(\alpha)$, and is independent of the
generator $\alpha$, since the ratio of any two generators is
a unit and has absolute norm $1$.

\begin{fact}\label{dets}
  With\/ $\N$ as above and $V\!$, $\phi$ as in the paragraph preceding
  Eq.~$\eqref{det}$, we have the identity\/ $\lvert\det^{}_{\QQ}
  (\phi)\rvert = \N \bigl(\det^{}_K (\phi) \bigr)$.
  In particular, $[I:J]=\N \bigl([I:J]^{}_K \bigr)$ for any
  ideals $I,J$ of $\HH (K)$ with $J\subseteq I$.
\end{fact}

\begin{proof}
  This is a standard result that can be proved by embedding $V$ in
  $V_{\RR}^r \oplus V_{\mathbb C}^s\simeq\mathbb R^{nd}$, where
  $V_{\mathbb R}$ is the real vector space $V\otimes_K\mathbb R$ (the
  \emph{realification} of $V$), $V_{\mathbb C}$ is the complex vector
  space $V\otimes_K\mathbb C$ (the \emph{complexification} of $V$),
  $r$ is the number of homomorphisms of $K$ into $\mathbb R$, $s$ is
  the number of complex conjugate pairs of homomorphisms of $K$ into
  $\mathbb C$, $d=r+2s$ is the degree of $K$ over $\QQ$ and $n$ is the
  dimension of $V$.  Then, $\phi$ can be extended to an $\RR$-linear
  map of $\mathbb R^{nd}$ into itself and the determinant of the
  extended linear map is the norm of the determinant of $\phi$.  The
  $1$-dimensional case of this construction is given in
  \cite[Ch.~II.3.1]{BS} and the general case in \cite[Sec.~4]{P}.
\end{proof}

This gives the following useful consequences.

\begin{prop}\label{N}
  Let $L'\subseteq L$ be arbitrary orders in\/ $\HH(K)$ and
  $q\in L^\bt$. Then, one has
\begin{itemize}
\item[\rm(i)]  $\quad[L:L'] \, = \, [\I(L):\I(L')];$
\smallskip
\item[\rm(ii)] $\quad[L:qL] \, = \, [L:Lq]
   \, = \, \N\bigl(\lvert q\rvert^4\bigr);$
\smallskip
\item[\rm(iii)]$\quad[\I(L):\I(qL)] \, = \, [\I(L):\I(Lq)]
   \, = \, \N\bigl(\lvert q\rvert^2\bigr),$
  provided that $q$ is $L$-primitive.
\end{itemize}
\end{prop}

\begin{proof}
Part(i) is immediate from Facts~\ref{im-index} and \ref{dets}.
Part~(ii) results from combining Fact~\ref{dets} with Lemma~\ref{[L:qL]},
and Part~(iii) from combining it with Lemma~\ref{[ImL:ImqL]}.
\end{proof}

\begin{remark}
Part (ii) of Proposition~\ref{N} clearly displays the correct
scaling behaviour of the index when $q$ is replaced by $\alpha q$
with $\alpha\in\oo$. Taking the scaling behaviour into account in
Part (iii), it is almost immediate that its extension to general
$q\in L^\bt$ reads
\begin{equation} \label{with-cont}
   [\I (L) : \I (qL)] \; = \; [\I (L) : \I (Lq)]
   \; = \;  \N\bigl( \mathrm{cont}^{}_{L} (q)\bigr)
   \N\bigl(\lvert q\rvert^2\bigr).
\end{equation}
It is possible \cite{Al} to derive generalizations of \eqref{with-cont}
to central simple algebras of dimension $n^2$ over fields with a
Dedekind ring $\oo$, with $q\mapsto\I(q)$ replaced by the mapping
$q\mapsto q-\frac{1}{n}\mathrm{tr} (q)$. Our formula would then
result from taking $n=2$ and considering only principal ideals.
\end{remark}

Let us now look at the coincidence problem for the rank $3$
$\oo$-module $\gG:=\I (L)$, where
we want to classify the intersections $\gG\cap R\gG$ with $R\in
\mathrm{SOC} (\gG)$ according to their indices. From
Proposition~\ref{Cayley-para}, we already know that, for finite
index, we may restrict to $R=R_q$ with $q\in L^\bt$, and with the
identity \eqref{matrix1} we have
\begin{equation} \label{intersection-one}
   \gG\cap R_q \gG \; = \; \I (L) \cap R_q \I (L)
   \; = \; \I (L) \cap \I (qL q^{-1}).
\end{equation}

\begin{lemma}\label{intersect}
  If $L$ is an order of\/ $\HH(K)$ and $q\in L^\bt$, one has
  the relation
\[
   \I(L)\cap \I(qLq^{-1}) \; = \; \I(L\cap qLq^{-1}).
\]
\end{lemma}

\begin{proof}
  The inclusion $\supseteq$ is immediate, so we have to prove the
  converse.  Consider an element $u^t\in\I(L)\cap \I(qLq^{-1})$, where
  we write $u^t$ for the transpose of $u$ without its first
  coordinate.  As $\I(qLq^{-1}) = R_q \I(L)$, there is some
  $v^t\in\I(L)$ with $u^t = R_q v^t$, and we can choose $x=(\alpha,u)$
  and $y=(\beta,v)$ in $L$, with $\alpha,\beta\in K$.  Assume for a
  moment that $\alpha-\beta\in\oo$. Then, $z=\alpha-\beta +y
  =(\alpha,v) \in L$ with $\I(z)=\I(y)$. Moreover, one has
\[
   q z q^{-1} \, = \, q (\alpha,v) q^{-1}
   \, = \, (\alpha,u) \, = \, x \, \in \, L
\]
so that $qzq^{-1} \in L\cap qLq^{-1}$ with $\I(qzq^{-1})=u^t$.

It remains to show that $\alpha-\beta\in\oo$. Observe that,
since $x$ and $y$ are integral,
\begin{align*}
\tr(x)&=2\alpha:=\gamma\in\oo,\\
\tr(y)&=2\beta:=\delta\in\oo, \mbox{ and}\\
\nr(x)-\nr(y)&=\alpha^2-\beta^2=(\gamma^2-\delta^2)/4\in\oo,
\end{align*}
the last line because $\nr(y)=\nr(qyq^{-1})=\nr(\beta,u)$.  We now
have $4|(\gamma^2-\delta^2)$ with $\gamma,\delta\in\oo$, and the proof
will be complete if we can deduce from this that $2|(\gamma-\delta)$.
Let $\pi$ be a prime factor of 2 and suppose that $\pi^e\parallel2$.
Now, $\pi^e\nmid(\gamma-\delta)$ would imply
$\pi^e\nmid(\gamma+\delta)$, since $\gamma+\delta\equiv\gamma-\delta$
(mod 2), giving $\pi^{2e}\nmid(\gamma^2-\delta^2)$ and contradicting
the fact that $\gamma^2-\delta^2$ is divisible by 4.  Hence every
prime factor of 2 divides $\gamma-\delta$ to at least as high a power
as it divides 2 and we have established that $2|(\gamma-\delta)$ and
thus $\alpha-\beta\in\oo$.
\end{proof}

The set of $L$-primitive elements of $L$ is still too big for our
purposes. Recall that any $q\in L^\bt$ with the property that $qL$ is a
two-sided ideal also satisfies $qLq^{-1}=L$, and is thus not of
interest for the coincidence problem (at least as far as counting the
CSMs is concerned).  In fact, the elements of $L$ that generate
two-sided ideals give, via the parametrization $R$, a multiple cover
of the rotation symmetry group of the module $\I(L)$, the multiplicity
being due to units of $\oo$.  Consider a general element $x\in L^\bt$
and assume that it is of the form $x=aqb$ with $a,b,q\in L^\bt$ and
$qL=Lq$, so that $q$ generates a 2-sided ideal. Then, there are
elements $a',b' \in L^\bt$ such that $x=qa' b = ab' q$, and $q$ is then
both a left and a right divisor of $x$. We can thus simply call $q$ a
divisor of $x$ in this case. Now, $x$ is called \emph{$L$-reduced} if
no divisor $q$ of $x$ other than the units of $L$ satisfies $qL=Lq$.
Again, we shall omit the prefix $L$ when it is clear from the context.
Recalling Eq.~\eqref{matrix1}, and thus the identity
\[
     R_q \I (L) \; = \; \I (qL q^{-1})\, ,
\]
the following consequence of Proposition~\ref{Cayley-para} and
Lemma~\ref{intersect} is immediate.

\begin{coro} \label{reduced-para}
  Let $\gG=\I (L)$ and consider a CSM of $\gG$, obtained from a rotation
  matrix $M\in\mathrm{SOC} (\gG) = \mathrm{SO} (3,K)$.  Then, among the\/
  $\mathrm{SOC} (\gG)$-matrices that result in the same CSM, there is at
  least one of the form $R=R(q)$ with $q$ an $L$-reduced element of $L$.
  \qed
\end{coro}

Consider now a maximal order $\OO$.  Our next aim is to relate the
conjugate orders $q\OO q^{-1}$, which are all of the same type as
$\OO$, to ideals $q\OO$ of $\OO$.  For this, we shall need to know
something about the irreducible factors of elements of $\OO$, and the
following fact is sufficient for our purpose.

\begin{fact}\label{factors}
  Let\/ $\HH(K)$ have class number $1$, and let $\OO$ be a maximal
  order in\/ $\HH (K)$. If $q\in\OO$ is primitive and $\pi\in\oo$ is a
  prime factor of\/ $|q|^2$, then $q$ can be factorized as $q=rp$,
  with $r,p\in\OO$ and\/ $|p|^2$ an associate of $\pi$.
\end{fact}

\begin{proof}
We translate the proof of \cite[Ch.~5, Thm.~2]{CSm} to this
more general setting.

Since $\HH(K)$ has class number $1$, the left ideal $\OO q+\OO\pi$ has
the form $\OO p$ for some $p\in \OO$ (the \emph{greatest common right
  divisor} of $q$ and $\pi$, determined up to multiplication by a
unit).  So $q=rp$ and $\pi=sp$, for some $r,s\in \OO$, and
$\pi^2=\nr(\pi)=|s|^2|p|^2$.  Since $\pi$ is prime, $|p|^2$ must be an
associate of 1, $\pi$ or $\pi^2$.  However, $p$ is not a unit because
$p=tq+u\pi$, for some $t,u\in \OO$, so
$\nr(p)=\nr(t)\nr(q)+\nr(u)\pi^2+\pi\tr(tq\bar u)$ is divisible by
$\pi$.  Also, if $|p|^2$ were an associate of $\pi^2$, then $s$ would
be a unit, making $\pi$ an associate of $p$ and therefore a divisor of
$q$, which is impossible since $q$ is primitive.  Consequently,
$|p|^2$ is an associate of $\pi$.
\end{proof}

We now embark on finding a correspondence between $q\OO$ and 
$\I (q\OO q^{-1})$ when $q$ is reduced.

\begin{lemma}\label{tricky}
  Let\/ $\HH(K)$ have class number $1$ and let $\OO$ be a maximal order.
  If $q,r\in\OO$ are primitive, one has
\[ \I(q\OO)\subseteq\I(r\OO)\;\Longrightarrow \; q\OO\subseteq r\OO
   \quad \mbox{and} \quad
   q\OO\subseteq r\OO \;  \Longrightarrow \;
  \OO\cap q\OO q^{-1}\subseteq \OO\cap r\OO r^{-1} .
\]
\end{lemma}

\begin{proof}
Since $\HH(K)$ has class number $1$, the right ideal $q\OO+r\OO$ has
the form $d\OO$ for some $d\in\OO$ (the \emph{greatest common left
divisor} of $q$ and $r$).  When $\I(q\OO)\subseteq\I(r\OO)$, we have
\[
   \I(d\OO)\; = \; \I(q\OO+r\OO)\; = \; 
   \I(q\OO)+\I(r\OO) \; = \; \I(r\OO),
\]
so, by Lemma~\ref{[ImL:ImqL]} and recalling $\lvert q\rvert^2 = \nr(q)$,
\[
   \nr(d)\, \oo \; = \; [\I(\OO):\I(d\OO)]^{}_K
   \; = \; [\I(\OO):\I(r\OO)]^{}_K \; = \; \nr(r)\, \oo,
\]
since $d$ (being a left divisor of $r$) is primitive.
Hence $\nr(r/d)$ is a unit in $\oo$, so $r/d$ is a unit in $\OO$ and
$r\OO=d\OO\supseteq q\OO$, establishing the first implication.

Now assume that $q\OO\subseteq r\OO$.  Then $q=rp$ for some $p\in\OO$.
Given any $x\in\OO\cap q\OO q^{-1}$, put $y=r^{-1}xr$ and define
\[
   \lden(y) \; := \; \{a\in\OO\mid ay\in\OO\}
\]
(the \emph{left denominator} of $y$).  This is a left ideal of $\OO$
and has the form $\OO d$ for some $d\in\OO$.  Clearly, $\lden (y)$
contains $r$, but it also contains $\bar{p}$, as $y=pzp^{-1}$ for some
$z\in\OO$.  Consequently, $r=sd$ and $\bar{p}=td$ with $s,t\in\OO$, and
hence $q=\lvert d\rvert^2 s\bar{t}$.  Since $q$ is primitive, $\lvert
d\rvert^2$ is a unit of $\oo$ and hence $d$ is a unit of $\OO$.  Thus
$\lden(y)=\OO$, so $y\in\OO$ and $x=ryr^{-1}\in\OO\cap r\OO r^{-1}$.
This establishes the second implication.
\end{proof}

\begin{lemma}\label{ideal-map}
  If\/ $\HH(K)$ has class number $1$, $\OO$ is a maximal order, 
  and $q\in\OO$ is reduced, then $\I(q\OO)=\I(\OO\cap q\OO q^{-1})$.
\end{lemma}

We note that, except when $q$ is a unit, it is only
the projections into $\HH_{-}$ that are equal:\ $q\OO$
and $\OO\cap q\OO q^{-1}$ themselves are distinct (because
$1\in q\OO$ only for $q$ a unit).

\begin{proof}
  We use induction on the absolute norm $\N\bigl(|q|^2\bigr)$.  The result is
  certainly true when $\N\bigl(|q|^2\bigr)=1$, since then $|q|^2$ is a unit in
  $\oo$ so $q$ is a unit in $\OO$ and $q\OO=\OO=\OO\cap q\OO q^{-1}$.  When
  $\N\bigl(|q|^2\bigr)>1$, $|q|^2$ is divisible by some prime $\pi$ of $\oo$,
  so $q=rp$ for some $r,p\in\OO$ with $|p|^2$ an
  associate of $\pi$, by Fact~\ref{factors}.  Moreover, $r$ and $p$ are
  reduced, because $q$ is. Now
\begin{equation}\label{sandwich}
   q\OO \, \subseteq \, \OO\cap q\OO q^{-1}
   \, \subseteq \, \OO\cap r\OO r^{-1},
\end{equation}
where the first inclusion is straightforward and the second is a
consequence of Lemma~\ref{tricky}.  Also $p\OO p^{-1}\neq \OO$, since
$p$ is reduced but not a unit, and  hence $q\OO q^{-1} \neq r\OO r^{-1}$.
But $\OO$, $q\OO q^{-1}$ and $r\OO r^{-1}$ are maximal orders,
so by Lemma~\ref{appendix} of the Appendix
$\OO\cap q\OO q^{-1}\ne\OO\cap r\OO r^{-1}$ whence
$[\OO\cap r\OO r^{-1}:\OO\cap q\OO q^{-1}]^{}_K \neq \oo$.
On projecting into $\HH_-$, (\ref{sandwich}) gives
\[
   \I(q\OO) \, \subseteq \, \I(\OO\cap q\OO q^{-1})
   \, \subseteq \, \I(\OO\cap r\OO r^{-1}) \, = \, \I(r\OO),
\]
where the equation on the right is an application of the induction 
hypothesis since $r$ is reduced and
\[
   \N\bigl(|r|^2\bigr) \, = \, 
   \N\bigl(|q|^2\bigr)/ \N\bigl(|p|^2\bigr)
   \, < \, \N\bigl(|q|^2\bigr).
\]
By Lemma~\ref{[ImL:ImqL]} and the multiplicativity
of \mbox{$K$-indices}, we have
\[
   [\I(r\OO):\I(q\OO)]^{}_K \, = \,
   \frac{[\I (\OO) :\I(q\OO)]^{}_K}{[\I (\OO) :\I(r\OO)]^{}_K}
   \, = \, \frac{ |q|^2}{|r|^2}\, \oo \, = \, \pi \oo \, ,
\]
a prime ideal of $\oo$, but Fact~\ref{im-index} implies
\[
   [\I(\OO\cap r\OO r^{-1}):\I(\OO\cap q\OO q^{-1})]^{}_K
  \; = \; [\OO\cap r\OO r^{-1}:\OO\cap q\OO q^{-1}]^{}_K 
   \; \neq \; \oo \, .
\] 
Being a divisor of the prime ideal $\pi\oo$, this index must
therefore be $\pi\oo$ itself and, again by the multiplicativity
of \mbox{$K$-indices}, $[\I(\OO\cap q\OO q^{-1}):\I(q\OO)]^{}_K=\oo$.
Hence $\I(q\OO)=\I(\OO\cap q\OO q^{-1})$, completing the induction.
\end{proof}

Summing up:

\begin{prop}\label{inclusions}
  If $K$ and\/ $\HH(K)$ both have class number $1$, if $\OO$ is a 
  maximal order of\/ $\HH (K)$, and if $q,r\in\OO$
  are \mbox{$\OO$-reduced}, the following inclusions are equivalent:
\begin{itemize}
\item[(i)]$\I(q\OO)\subseteq\I(r\OO);$
\item[(ii)]$q\OO\subseteq r\OO;$
\item[(iii)]$\OO\cap q\OO q^{-1}\subseteq\OO\cap r\OO r^{-1};$
\item[(iv)]$\I(\OO\cap q\OO q^{-1})\subseteq\I(\OO\cap r\OO r^{-1})$.
\end{itemize}
\end{prop}

\begin{proof}
Lemma~\ref{tricky} gives $\mathrm{(i)\Rightarrow(ii)\Rightarrow(iii)}$,
the claim $\mathrm{(iii)\Rightarrow(iv)}$ follows by projection into 
$\HH_-$, and $\mathrm{(iv)\iff(i)}$ by Lemma~\ref{ideal-map}.
\end{proof}

As a matter of independent interest, we list the following consequences of
Lemma~\ref{ideal-map}.

\begin{prop} \label{prop5}
  Let\/ $\HH(K)$ be the standard quaternion algebra
  over the real algebraic number field $K$, and suppose
  that both $K$ and\/ $\HH(K)$ have class number $1$.
  Let $\OO$ be a maximal order of\/ $\HH(K)$
  and $q\in\OO$ an $\OO$-reduced element. Then, one has the
  following relations:
\begin{itemize}
\item[\rm(i)]   $\quad\I(\OO\cap q\OO q^{-1}) \, = \, \I(q\OO)
     \, = \, \I(\OO \bar q);$
\smallskip
\item[\rm(ii)]  $\quad\OO\cap q\OO q^{-1} \, = \, \oo+q\OO+\OO\bar q
     \, = \, \oo+q\OO \, = \, \oo+\OO\bar q;$
\smallskip
\item[\rm(iii)] $\quad[\OO:\OO\cap q\OO q^{-1}] \, = \,
     [\I(\OO):\I(\OO\cap q\OO q^{-1})] \, = \,
     \N\bigl(\lvert q\rvert^2\bigr);$
\smallskip
\item[\rm(iv)] $\quad[\I(\OO):\I(\OO)\cap \I(q\OO q^{-1})] \, = \,
     \N\bigl(\lvert q\rvert^2\bigr);$
\smallskip
\item[\rm(v)]  $\quad[\OO\cap q\OO q^{-1}:q\OO] \, = \,
     [\oo+q\OO:q\OO] \, = \, [\oo+\OO\bar q:q\OO]
     \, = \, \N\bigl(\lvert q\rvert^2\bigr)$.
\end{itemize}
\end{prop}

\begin{proof}
Part~(i) follows from Lemma~\ref{ideal-map} together with
the observation that $\I(qx)=\I(-\bar x\bar q)$ when $x\in\OO$.
Part~(ii) follows from the inclusions
\[
   \oo+q\OO \;\subseteq\; \oo+q\OO+\OO\bar q
   \;\subseteq\; \OO\cap q\OO q^{-1}
\]
and the fact that, by Propositions~\ref{indices} and Lemma~\ref{ideal-map},
\[
    [\OO\cap q\OO q^{-1}:\oo+q\OO]^{}_K \; = \;
    [\oo:\oo]^{}_K \, [\I(\OO\cap q\OO q^{-1}):\I(q\OO)]^{}_K
    \; = \; \oo .
\]
For Part~(iii), we note that
\[
     [\OO:\OO\cap q\OO q^{-1}]^{}_K \; = \;
     [\oo:\oo]^{}_K \, [\I(\OO):\I(\OO\cap q\OO q^{-1})]^{}_K,
\]
by Proposition~\ref{indices}, and use Lemmas~\ref{ideal-map} and 
\ref{[ImL:ImqL]}, and Fact~\ref{dets}.  Part~(iv) follows
from Part~(iii) and Lemma~\ref{intersect}. Finally, Part~(v)
is a consequence of Parts~(ii) and (iii), Lemma~\ref{[L:qL]},
the multiplicativity of $K$-indices and Fact~\ref{dets}.
\end{proof}

\section{General results} \label{general}

We are finally in a position to state our main result on the
connection between coincidence site modules and one-sided ideals.

\begin{theorem} \label{main-result}
  Let\/ $\HH(K)$ be the standard quaternion algebra over
  the real algebraic number field $K$, and suppose that both 
  $K$ and\/ $\HH (K)$ have class number $1$.
  Let\/ $\OO$ be a maximal order of\/ $\HH(K)$.
  Then, the CSMs of\/ $\gG := \I (\OO)$ are in one-to-one
  correspondence with the one-sided ideals $q\OO$ that are
  generated by the \mbox{$\OO$-reduced} elements\/ $q\in\OO$.
  The corresponding coincidence indices are given by\/
  $\N(\vert q\rvert^2)$, so that these right ideals are to be
  counted according to the square root of their ideal index
  from Proposition~$\ref{N}$ {\rm (ii)}.
\end{theorem}

\begin{proof}
Parts (ii) and (iv) of
Proposition~\ref{inclusions} provide the connection between
the CSMs and the right ideals generated by \mbox{$\OO$-reduced}
elements of $\OO$, and Proposition~\ref{N} (iii) gives the indices
of the CSMs.
\end{proof}

This implies that we may count the CSMs by counting one-sided ideals
in $\OO$, which can be done by counting \emph{all} ideals of a given
index in $\OO$ and subtracting off the two-sided ones. While the
former is achieved by the zeta function $\zeta^{}_{\OO} (s)$, the latter
is done by the zeta function of the base field, via $\zeta^{}_{K} (4s)$,
up to finitely many correction factors from ramified primes. The
factor $4$ in the argument of $\zeta^{}_{K}$ stems from the index
formula of Lemma~\ref{[L:qL]}. Let us introduce a Dirichlet series
generating function for one-sided ideals that are generated by
\mbox{$\OO$-reduced} elements as $\zeta^{\rm red}_{\OO}(s)$, and
denote the Dirichlet series generating function for the two-sided
ideals of $\OO$ by $\zeta^{}_{\OO.\OO} (s)$, see below for examples.
\begin{lemma} \label{zeta-reduced}
  If $\OO$ is a maximal order in the quaternion algebra\/
  $\HH(K)$  of class number\/ $1$, one has the identity
  $\; \zeta^{}_{\OO} (s) =  \zeta^{\rm red}_{\OO} (s)
   \, \zeta^{}_{\OO.\OO} (s) $.
\end{lemma}

\begin{proof}
Due to the assumption on the class number, each right ideal is
principal, and thus of the form $q\OO$ for some $q\in\OO$.  If $q$ has
any divisor that generates a two-sided ideal, we can pull that factor
out to the right to obtain a unique factorization of $q\OO$ into a
reduced and a two-sided ideal. Consequently, sorted by their index,
we can either count all right ideals together (giving the term on the
left), or reduced right ideals and their possible right multiplication
by two-sided ideals separately (giving the product of Dirichlet series
on the right).
\end{proof}

Now, in view of Theorem~\ref{main-result}, counting CSMs is the same
as counting reduced right ideals, but with respect to the square root
of the ideal index. This results in
\begin{theorem} \label{csm-diri}
  Under the assumptions of Theorem~$\ref{main-result}$, the Dirichlet
  series generating function\/ $\varPhi (s)$ for the number\/ $f(m)$
  of CSMs of index\/ $m$ of the\/ $\oo$-module\/ $\gG=\I (\OO) \subset
  \RR^3$ is given by
\[
    \varPhi (s) \; = \; \sum_{m\ge 1} \frac{f(m)}{m^s} \; = \;
    \zeta^{\rm red}_{\OO} (s/2) \; = \;
    \frac{\zeta^{}_{\OO} (s/2)}{\zeta^{}_{\OO.\OO} (s/2)}
    \; = \;  E(s)\, \frac{\zeta^{}_{K}(s)\,
    \zeta^{}_{K}(s-1)}{\zeta^{}_{K} (2s)}\, ,
\]
   where $\zeta^{}_{K} (s)$ is the Dedekind zeta function of the field
   $K$ and $E(s)$ is either $1$ or an additional analytic factor $($with
   finitely many terms\/$)$ that takes care of the extra contributions
   from ramified primes.
   
   In particular, the arithmetic function\/ $f(m)$ is multiplicative, and
   the elementary coincidence spectrum of\/ $\gG$ is\/
   $\varSigma(\mathrm{OC} (\gG))=\varSigma(\mathrm{SOC} (\gG))= \{
   m\in\NN \mid f(m)\neq 0\}$.
\end{theorem}

\begin{proof}
The first claim follows from Theorem~\ref{main-result} together with
Lemma~\ref{zeta-reduced}. Since both zeta functions involved,
$\zeta^{}_{\OO} (s)$ and $\zeta^{}_{\OO.\OO} (s)$, possess an Euler
product decomposition, this applies to $\varPhi(s)$ as well, with
their detailed structure following from \cite{V} and \cite{R}.  As a
consequence, $f(m)$ must be a multiplicative arithmetic function.  The
claim about the spectrum is obvious.
\end{proof}

\begin{coro} \label{asymp-general}
  The asymptotic growth of the multiplicative arithmetic function\/ $f(m)$
  of Theorem~$\ref{csm-diri}$ is given through its summatory function via
\[
   F(x) \; := \; \sum_{m\le x} f(m) \; \sim \;  \varrho\, \frac{x^2}{2} ,
   \quad \mbox{as } x\to\infty \ts , \quad
   \mbox{with }\,
   \varrho = \mathrm{res}^{}_{s=2}\, \varPhi (s)\ts .
\]   
   In particular, the average size of $f(m)$ is $\varrho\ts m$.
\end{coro}

\begin{proof}
Consider the function $\varPhi (s)$ from Theorem~\ref{csm-diri}. 
By the general structure of the zeta functions involved and the fact
that $E(s)$ is an entire function, $\varPhi(s)$ is
analytic in the open half-plane $\{ s=\sigma + i t \mid \sigma > 2\}$,
and has a first order pole at $s=2$, but no other singularity on the
line $\{ \sigma = 2\}$. In fact, $\varPhi (s)$ is analytic everywhere
on this line except at $s=2$.

Delange's theorem, compare \cite[Appendix]{BM} for a formulation tailored
to this situation, now gives $F(x) \sim \varrho\ts x^2/2$, as $x\to\infty$,
with $\varrho$ as claimed.
\end{proof}

\begin{remark}
The residue in Corollary~\ref{asymp-general} can be calculated as
\begin{equation} \label{residue}
   \mathrm{res}^{}_{s=2}\, \varPhi (s)
   \; = \; \frac{E(2)\,\zeta^{}_{K} (2)}{\zeta^{}_{K} (4)}\,
   \mathrm{res}^{}_{s=1}\, \zeta^{}_{K} (s) \ts .
\end{equation}
and $\mathrm{res}^{}_{s=1}\, \zeta^{}_{K} (s)$
can be expressed in terms of values of $L$-series.
\end{remark}

\section{Specific results}
\label{specific}

There are three explicit cases of particular interest for applications
in crystallography, where the ring of integers $\oo$ is a principal
ideal domain (compare Table~1 on p.~422 of \cite{BS}) and the
quaternion algebra $\HH(K)$ has class number $1$ (see \cite{St} and
the first table on p.~156 of \cite{V} for further examples).  The
basic information is summarized in Table~\ref{tab}. The corresponding
zeta functions are the Dirichlet series generating functions for the
right ideals,
\[
    \zeta^{}_{\OO} (s) \; = \; \sum_{0\neq\mathfrak{A}\subset\OO}
    \frac{1}{[\OO:\mathfrak{A}]^s} \; = \; \sum_{n\ge 1}
    \frac{f^{}_{\OO}(n)}{n^s} \, ,
\]
where $\mathfrak{A}$ runs through the non-zero right ideals of
$\OO$ and $f^{}_{\OO} (m)$ is the number
of such ideals of index $m$. Due to the class number being $1$, all
ideals are principal, so we have $\mathfrak{A}={a}\OO$ for some
${a}\in\OO^\bt$, and then $[\OO:{a}\OO]=\N\bigl(\lvert{a}\rvert^4\bigr)$
by Proposition~\ref{N}~(ii). The zeta function for the left ideals is
the same, as left and right ideals are in a natural one-to-one relation
that preserves the norm. Similarly, $\zeta^{}_{\OO . \OO} (s)$ is
defined as the Dirichlet series that runs over all non-zero two-sided
ideals of $\OO$, so that  $\zeta^{}_{\OO.\OO} (s) = \sum_{n\ge 1}
\frac{f_{\OO.\OO} (n)}{n^s}$.

{\large \begin{table}
\begin{tabular}{c|ccc}
$\quad K \quad$   & $\quad \QQ \quad$ &
$\QQ(\sqrt{5}\,)$ & $\quad\QQ(\sqrt{2}\,)\quad$
\rule[-2ex]{0ex}{5ex} \\ \hline
$\OO$ & $\JJ$ & $\II$              & $\KK$
\rule[-2ex]{0ex}{5ex} \\
$\oo$ & $\ZZ$ & $\ZZ[\tau]$        & $\ZZ[\sqrt{2}\,]$
\rule[-2ex]{0ex}{4ex}
\end{tabular}
\bigskip
\caption{Data for three quaternion algebras $\HH(K)$ of class
number $1$, with maximal order $\OO$. Here, $\tau=(\sqrt{5}+1)/2$
is the golden ratio, and the ring of integers $\oo$ is a PID in all
three cases (see text for details). \label{tab}}
\end{table}}

In all three cases of Table~\ref{tab}, these zeta functions can be expressed
in terms of the Dedekind zeta functions of $K$ as follows, compare \cite{V}:
\begin{equation} \label{zeta1}
   \zeta^{}_{\OO} (s) \; = \;
   \zeta^{}_{K} (2s)\, \zeta^{}_{K} (2s-1) \cdot
   \begin{cases} (1-2^{1-2s}), & \text{if $K=\QQ$}, \\
   1, & \text{otherwise}, \end{cases}
\end{equation}
where the extra factor for $K=\QQ$ results from the rational prime $2$
being ramified in this case. The Dedekind zeta functions entering
are given by
\begin{eqnarray} \label{zeta3}
   \zeta^{}_{\QQ} (s) & = &
   \zeta (s) \; = \; \prod_{p\in\mathcal{P}}
   \frac{1}{1-p^{-s}}\, , \nonumber \\
   \zeta^{}_{\QQ(\sqrt{5}\,)} (s) & = & \frac{1}{1-5^{-s}}
   \prod_{p\equiv\pm 1 \; (5)} \frac{1}{(1-p^{-s})^2}
   \prod_{p\equiv\pm 2 \; (5)} \frac{1}{1-p^{-2s}}\, , \\
   \zeta^{}_{\QQ(\sqrt{2}\,)} (s) & = & \frac{1}{1-2^{-s}}
   \prod_{p\equiv\pm 1 \; (8)} \frac{1}{(1-p^{-s})^2}
   \prod_{p\equiv\pm 3 \; (8)} \frac{1}{1-p^{-2s}}\, , \nonumber
\end{eqnarray}
see \cite{BM} for further details in this context, including explicit
terms and asymptotic properties.

We shall also need the Dirichlet series generating functions for the
two-sided ideals of $\OO$. By the relation between prime ideals in
$\OO$ and $\oo$ according to \cite[Thm.~22.4]{R}, it can be expressed
by the zeta function of the base field $K$, possibly up to finitely
many correction factors for primes that ramify in the extension from
$K$ to $\HH(K)$. With the examples of Table~\ref{tab}, this happens
only for $\JJ$, where $(1+\ii)\JJ = \JJ (1+\ii) = \JJ (1+\ii) \JJ$ is the
unique two-sided prime ideal over the rational prime $2$. We know from
Proposition~\ref{N}~(ii) that $(1+\ii)\JJ$ has index $4$ in $\JJ$, to be
compared with index $16$ for the ideal $2\JJ$. In all other cases, if
$\mathfrak{P}$ is the unique two-sided prime ideal in $\OO$ that lies
over the prime ideal $\mathfrak{p}\in\oo$, compare
\cite[Thm.~22.3]{R}, one has
\[
     [\OO : \mathfrak{P} ] \; = \; [\oo : \mathfrak{p} ]^4 .
\]
Consequently, we obtain for our three examples of Table~\ref{tab}:
\begin{equation} \label{zeta-two-sided}
   \zeta^{}_{\OO.\OO} (s) \; = \;
   \zeta^{}_{K} (4s) \cdot
   \begin{cases} (1+4^{-s}), & \text{if $K=\QQ$}, \\
   1, & \text{otherwise}. \end{cases}
\end{equation}
With these general preparations, we can now apply the results to the
three cases of Table~\ref{tab}.

\subsection{The Hurwitz ring $\JJ$} \label{Hurwitz-results}
Here, we consider $\HH(\QQ)$, with $\JJ$ as defined in Eq.~\eqref{J-def}.
This is the unique maximal order that contains $\LL\simeq \ZZ^4$ of
\eqref{lattice}, with $\oo=\ZZ$ and $[\JJ:\LL]=2$. It is immediate that
$\I (\JJ) = \gG_{\rm bcc}$, so that we can re-derive the results from
Section~\ref{example}. The characterization of the \mbox{$\JJ$-reduced}
elements of $\JJ$ follows from the detailed analysis in \cite{H}:

\begin{fact} \label{J-reduced}
  An element\/ $q\in\JJ$ is \mbox{$\JJ$-reduced} if and only if\/
  $q$ is\/ $\JJ$-primitive and\/ $\vert q\rvert^2$ is odd. \qed
\end{fact}

Theorem~\ref{csm-diri}, applied to $\JJ$, reproduces all results of
Section~\ref{example}. In particular, we find

\begin{coro}
  The elementary coincidence spectrum of the lattice $\gG_{\rm bcc}=\I
  (\JJ)$ is the set of odd integers. This set is a monoid, generated
  by the odd rational primes.
  \qed
\end{coro}

In view of Fact~\ref{J-reduced}, one can also derive a simple formula
for the index of a rotation matrix $R\in\mathrm{SOC} (\gG_{\rm bcc})$.
If $R=R_q$ with $q\in\JJ$ being $\JJ$-reduced, one simply has
$\varSigma_{\rm bcc} (R_q) = \lvert q\rvert^2$. Otherwise, replace $q$ by
$q'=q/\mathrm{cont}^{}_{\JJ} (q)$, which is $\JJ$-primitive. If
$\lvert q'\rvert^2$ is odd, the quaternion $q'$ is already reduced.
If not, one has $q'=q'' (1+\ii)$ with $\lvert q''\rvert^2$ odd.
In this case, since $1+\ii$ generates a two-sided ideal, $R_{q''}$ defines
the same CSL as $R_{q'}$.  For a general $q\in\JJ$, we thus have
\begin{equation} \label{matrix-index1}
    \varSigma_{\rm bcc} (R_q) \; = \;
    \bigl\lvert q/\mathrm{cont}^{}_{\JJ} (q) \big\rvert^2 \cdot
    \begin{cases}  1, & \text{if $\lvert q'\rvert^2$ is odd}, \\
                              \tfrac{1}{2}, & \text{otherwise}.
     \end{cases}                         
\end{equation}
This formula is well known and has been derived before \cite{G,B} by 
different methods.

The dual lattice of $\gG_{\rm bcc}$ is the so-called face centred
cubic lattice, conveniently written as
\[
    \gG_{\rm fcc} \; = \; \{ (\ell,m,n)^t \in \ZZ^3\mid
    \ell + m + n \;\mbox{ even}\} \; = \;
    (1+\ii) \JJ \cap \HH_- \, .
\]
It is a sublattice of $\ZZ^3$ of index $2$. Clearly, due to
commensurability, all three lattices share the same commensurator
subgroup of $\mathrm{SO} (3,\RR)$, namely $\mathrm{SO} (3,\QQ)$.
More generally, we have
\begin{theorem} \label{cubic-index}
  Let\/ $R\in\mathrm{SO} (3,\QQ)$. Then, the three cubic lattices in\/
  $\RR^3$ share the same index formula, i.e., $\varSigma_{\rm bcc} (R) =
  \varSigma_{\rm p} (R) = \varSigma_{\rm fcc} (R)$.  In
  particular, the coincidence spectrum for all three lattices is the set of
  odd integers, as given in Eq.~$\eqref{spec1}$, and the Dirichlet series
  generating function for all three cubic lattices is\/ $\varPhi^{}_{\rm
    cub} (s)$ of Eqs.~$\eqref{bcc-diri}$ and $\eqref{J-diri}$.
\end{theorem}

\begin{proof}
By \cite[Thm.~2.2]{B}, two mutually dual lattices share the same
index formula, which then applies to $\gG_{\rm bcc}$ and $\gG_{\rm
  fcc}$. So, we have $\varSigma_{\rm bcc} (R) = \varSigma_{\rm fcc}
(R) =: \varSigma (R)$.

Moreover, we have the inclusion
\[
    \gG_{\rm fcc} \;\stackrel{2}{\subset}\; \ZZ^3
    \;\stackrel{2}{\subset}\; \gG_{\rm bcc}\, .
\]
By \cite[Lemma~2.6]{B}, this relation implies $\varSigma_{\rm p}
(R)\mid 2 \varSigma (R)$ and $\varSigma (R)\mid 2 \varSigma_{\rm p}
(R)$, which tells us that the coincidence index for $\ZZ^3$ could only
differ from $\varSigma (R)$ by a factor of $2$.

{}From our previous arguments, we know that $\varSigma (R)$ is always
odd, so that $\varSigma (R)\mid \varSigma_{\rm p} (R)$ On the other
hand, also the index $\varSigma_{\rm p} (R)$ is always odd, see
\cite{G,B}, which emerges from an independent argument directly on the
basis of the rotation matrix $R$. Consequently, we must have
$\varSigma_{\rm p} (R)\mid \varSigma (R)$, which gives $\varSigma_{\rm
  p} (R) = \varSigma (R)$.
\end{proof}

\begin{remark}
Observing that the lattice $\ZZ^3$ is isomorphic with $\JJ\cap\HH_{-}$,
the index relation $\varSigma_{\rm p} (R)= \varSigma (R)$ in the previous
proof can also be derived with an independent argument based upon
Fact~\ref{real-part-of-O} and an index formula for intersections with 
$\HH_-$, see \cite{BPnew} for more.
\end{remark}

Also of interest is the analysis of multiple intersections.
Unfortunately, the relation to the arithmetic of $\JJ$ becomes
considerably more involved. First explicit results are available
\cite{Z2,Z3}, including the observation that the total coincidence
spectrum is identical to the elementary one. For a given index,
new CSLs appear, but all possibilities are exhausted with
triple intersections.

\subsection{The icosian ring $\II$}
Here, we consider $\HH\bigl(\QQ(\sqrt{5}\,)\bigr)$. In this case,
there are two maximal orders that contain $\LL$ of \eqref{lattice},
with $\oo=\ZZ[\tau]$. One choice is
\[
   \II \; := \; \big\langle (1,0,0,0),(0,1,0,0),
   \tfrac{1}{2} (1,1,1,1),\tfrac{1}{2}
   (1\!-\!\tau,\tau,0,1) \big\rangle_{\ZZ[\tau]}\, ,
\]
where $\tau=(\sqrt{5}+1)/2$ is the golden ratio. The other choice is
$\II'$, where ${}'$ denotes the algebraic conjugation of
$\QQ(\sqrt{5}\,)$, defined by $\sqrt{5}\mapsto -\sqrt{5}$. Both have been
called the \emph{icosian ring}, see \cite{CMP,BM,CSl} for details on the
connection to the root system of type $H_4$. Our choice here matches
that of \cite{CMP,CSl}, while the roles of $\II$ and $\II'$ are
interchanged in comparison to \cite{B,BM}. It is not difficult to
check that $\LL$ of \eqref{lattice} is an index $16$ submodule both of
$\II$ and $\II'$. Moreover, one has the relation
\[
   2\ts ( \II + \II' ) \; \stackrel{4}{\subset} \;
   \LL \; \stackrel{4}{\subset} \; \II \cap \II'  ,
\]
where $\II \cap \II' = \JJ [\tau] := \JJ + \tau\JJ$ provides a
connection to the Hurwitz ring.  Let us mention that, in the icosian
case, $2\II$ is the unique prime ideal over $2\in\ZZ[\tau]$, and that
$(1+\ii)\II = \II' (1+\ii)$, so that $(1+\ii)\II$ is a one-sided, but not a
two-sided ideal. This also shows that $\II$ and $\II'$ are related by an
inner automorphism, and hence are of the same type, in line with
Fact~\ref{types-and-classes}.  Unlike in the case of the Hurwitz ring,
the rational prime $2$, which is also a prime of $\ZZ[\tau]$, does not
ramify. More generally, no prime of $\ZZ[\tau]$ ramifies in the
extension step to $\II$.  Consequently, one has
\begin{fact}
  An element $q\in\II$ is\/ $\II$-reduced if and only if 
  it is\/ $\II$-primitive.  \qed
\end{fact}

With this definition, the imaginary part satisfies $\I (\II)=
\frac{1}{2}\mathcal{M}_{\rm B}$, where
\[
   \begin{split}
   \mathcal{M}_{\rm B}  &\; = \; \big\{
   \mbox{$\sum_{i=1}^{3}$ } \alpha_i e_i
   \mid \alpha_i \in\ZZ[\tau], \mbox{ with }
   \tau^2\alpha^{}_{1} + \tau\alpha^{}_{2} + \alpha^{}_{3}
   \equiv 0 \; (\mbox{mod } 2) \big\} \\[2mm]
   & \; = \; \big\langle (2,0,0)^t, (1,1,1)^t,
   (\tau,0,1)^t \ts\big\rangle_{\ZZ[\tau]} \; \subset \; \RR^3
   \end{split}
\]
is the standard body centred icosahedral module of
quasi-crystallography, see \cite{B} and references therein.
Furthermore, $\II\cap\HH_{-}\simeq
\I (\II\cap\HH_{-}) =\tfrac{1}{2}\mathcal{M}_{\rm F}$, 
where
\[
   \begin{split}
   \mathcal{M}_{\rm F}  &\; = \; 
   \bigl\{
   \mbox{$\sum_{i=1}^{3}$ } \alpha_i e_i
   \mid \alpha_i \in\ZZ[\tau], \mbox{ with }
   \alpha^{}_{1} \equiv \tau\alpha^{}_{2} \equiv
   \tau^2\alpha^{}_{3} \; (\mbox{mod } 2) \bigr\} \\[2mm]
   & \; = \; \bigl\{
   \mbox{$\sum_{i=1}^{3}$ } \alpha_i e_i \in \mathcal{M}_{\rm B}
   \mid 
   \alpha^{}_{1} + \alpha^{}_{2} + \alpha^{}_{3} \equiv 0 
   \; (\mbox{mod } 2) 
   \bigr\} \\
   & \; = \; \bigl\langle (2,0,0)^t, (\tau \! + \! 1,\tau,1)^t,
   (0,0,2)^t \ts\bigr\rangle_{\ZZ[\tau]} 
   \; \stackrel{4}{\subset} \; \mathcal{M}_{\rm B}
   \; \subset \; \RR^3
   \end{split}
\]
is the standard face centred icosahedral module. Both
$\mathcal{M}_{\rm B}$ and $\mathcal{M}_{\rm F}$ are free
$\ZZ[\tau]$-modules of rank $3$, and free $\ZZ$-modules of rank
$6$. One quickly checks that
\[
   [ \I (\II) : \I (\LL) ] \; = \;
   [ \II : \LL ] \; = \; 16\, ,
\]
in line with Fact~\ref{im-index}. Note that $\I (\II)$ can be
obtained as the projection of the weight lattice $D^*_6$ into a
$3$-dimensional subspace that is invariant under the symmetry 
group of the icosahedron. This also relates $\mathcal{M}_{\rm B}$
to a lattice in $6$-space, which is complemented by the corresponding
relation between $\mathcal{M}_{\rm F}$ and the root lattice $D_6$,
see \cite{B} for details.

An application of Theorem~\ref{csm-diri} to $\II$ gives the generating
function
\begin{equation} \label{I-diri}
   \begin{split}
   \varPhi_{\rm ico} (s) &  \; = \; \frac{1+5^{-s}}{1-5^{1-s}}
   \prod_{p\equiv\pm 1 \; (5)} \left(\frac{1+p^{-s}}{1-p^{1-s}}\right)^2
   \prod_{p\equiv\pm 2 \; (5)} \frac{1+p^{-2s}}{1-p^{2(1-s)}}  \\
   & \; = \; 1 + \frac{5}{4^s} + \frac{6}{5^s} + \frac{10}{9^s} +
   \frac{24}{11^s} + \frac{20}{16^s} + \frac{40}{19^s} + \frac{30}{20^s}
   + \frac{30}{25^s} + \frac{60}{29^s} + \ldots
   \end{split}
\end{equation}

\begin{remark} \label{remark5}
  The generating function $\varPhi_{\rm ico} (s)$ applies to both
  $\mathcal{M}_{\rm B}$ and $\mathcal{M}_{\rm F}$, as well as to any image of
  either under a linear similarity.  The relation between $\mathcal{M}_{\rm
    B}$ and $\mathcal{M}_{\rm F}$ can be understood from the point of view of
  mutually dual embeddings as lattices into $6$-space mentioned before,
  compare \cite{B}.
\end{remark}

The asymptotic growth constant $\varrho$ of
Corollary~\ref{asymp-general} takes the form
\[
    \varrho \; = \; \mathrm{res}_{s=2} \,\varPhi_{\rm ico} (s) 
    \; = \; \frac{45\sqrt{5} \log (\tau)}{\pi^4}
    \; \simeq \; 0.497\, 089\, ,
\]
and $\varrho\ts m$ is the average size of the arithmetic function $f(m)$;
details of the calculation follow from \cite[Appendix]{BM}.

The index of a matrix $R\in\mathrm{SOC} (\mathcal{M}_{\rm B})$,
written as $R=R_q$ with $q\in\II^\bt$, is given by
\begin{equation}
   \varSigma (R_q) \; = \;
   \N \bigl( \lvert q / \mathrm{cont}^{}_{\II} (q) \rvert^2 \bigr).
\end{equation}
It is thus of the form $\N (k +\ell\tau)=k^2 + k\ell - \ell^2$ with
$k,\ell\in\ZZ$.  Since each prime $\pi\in\ZZ[\tau]$ has a
decomposition $\pi=p\bar{p}$ in $\II$, one sees that all indices of
this form indeed occur.
\begin{coro}
  The\/ $\ZZ[\tau]$-modules $\mathcal{M}_{\rm B}$ and 
  $\mathcal{M}_{\rm F}$ have the\/ {\rm SOC}-group\/
  $\mathrm{SO} (3,\QQ(\sqrt{5}\,))$ and share the
  elementary coincidence spectrum
\[
   \begin{split}
   \varSigma \bigl(\mathrm{SO} (3,\QQ(\sqrt{5}\,))\bigr)
   & \; = \; \{ m\in\NN \mid m= k^2 + k\ell - \ell^2 \, ,
   \mbox{ with } k,\ell\in\ZZ \} \\[2mm]
   & \; = \;
   \{ 1,4,5,9,11,16,19,20,25,\ldots \}\, .
   \end{split}
\]
   This spectrum is a monoid and consists of all natural numbers that
   have all primes congruent to\/ $\pm 2$ $(\mbox{mod } 5)$ in their
   factorizations occuring with an even power.  
\end{coro}

\begin{proof}
The claim on $\mathcal{M}_{\rm B}$ is clear by Theorem~\ref{csm-diri},
as is the characterization of the natural numbers representable by
the quadratic form $k^2 + k\ell - \ell^2$, compare \cite{HW}.

The SOC-group of $\mathcal{M}_{\rm F}$ is the same because the two
modules are commensurate. Also, the index formula is the same because
the modules $\mathcal{M}^{}_{F}$ and $\mathcal{M}^{}_{B}$ are dual to
one another via the embedding into $6$-space (as $D^{}_{6}$ and
$D^{*}_{6}$, respectively). Consequently, one can use the result from
\cite{B} that they share the index formula, thus giving equal spectra.
\end{proof}

A related object of interest is $\ZZ[\tau]^3 = \I (\LL)$, which is a
submodule of $\I (\II)$ of index $16$ and hence commensurate. Clearly,
$\mathrm{OC} \bigl( \I (\LL)\bigr) = \mathrm{OC} \bigl( \I
(\OO)\bigr)$, and they also share the elementary coincidence spectrum
\cite{B}, though a given index may arise from different sets of
rotations.  Consequently, one gets different Dirichlet series, too,
see \cite[Sec.~5.3]{B} for details.

\subsection{The octahedral (or cubian) ring $\KK$}
Here, we consider $\HH\bigl(\QQ(\sqrt{2}\,)\bigr)$, with maximal order
\[
   \KK \; := \; \big\langle 1,\tfrac{1+\ii}{\sqrt{2}},
   \tfrac{1+\jj}{\sqrt{2}},\tfrac{1+\ii+\jj+\kk}{2}
   \big\rangle_{\ZZ[\sqrt{2}\,]} \, ,
\]
see \cite{BM} for a more detailed description. Note that $\KK$
contains $\LL$ of \eqref{lattice}, with $\oo=\ZZ[\sqrt{2}\,]$, as an
index $16$ submodule and that $\KK$ is the only maximal order to
contain $\LL$.  In this case, the imaginary part is
\[
   \I (\KK) \; = \; \tfrac{1}{\sqrt{2}}\, \big\langle
   (1,0,0)^t, (0,1,0)^t, \tfrac{1}{\sqrt{2}} (1,1,1)^t
   \big\rangle_{\ZZ[\sqrt{2}\,]}\, ,
\]
so that $[\I (\KK) : \I (\LL)] = [ \KK : \LL] = 16$, in line with
Fact~\ref{im-index}. In the case of $\KK$, the rational prime $2$
ramifies in the extension from $\ZZ$ to $\ZZ[\sqrt{2}\,]$, but
$\sqrt{2}$ does not ramify in the following extension to $\KK$. As no
other prime ramifies, one finds
\begin{fact}
  An element $q\in\KK$ is\/ $\KK\!$-reduced
  if and only if it is\/ $\KK\!$-primitive.  \qed
\end{fact}

This time, an application of Theorem~\ref{csm-diri} gives
\begin{equation} \label{K-diri}
   \begin{split}
   \varPhi_{\rm oct} (s) &  \; = \; \frac{1+2^{-s}}{1-2^{1-s}}
   \prod_{p\equiv\pm 1 \; (8)} \left(\frac{1+p^{-s}}{1-p^{1-s}}\right)^2
   \prod_{p\equiv\pm 3 \; (8)} \frac{1+p^{-2s}}{1-p^{2(1-s)}}  \\
   & \; = \; 1 + \frac{3}{2^s} + \frac{6}{4^s} + \frac{16}{7^s} +
   \frac{12}{8^s} + \frac{10}{9^s} + \frac{48}{14^s} + \frac{24}{16^s}
   + \frac{36}{17^s} + \frac{30}{18^s} + \ldots
   \end{split}
\end{equation}
Here, the asymptotic growth constant $\varrho$ of
Corollary~\ref{asymp-general} reads
\[
    \varrho \; = \; \mathrm{res}_{s=2}\, \varPhi_{\rm oct} (s)
    \; = \; \frac{720 \sqrt{2} \log (1+\sqrt{2}\,)}{11\,\pi^4}
    \; \simeq \; 0.837\, 559 \, ,
\]
with $\varrho\ts m$ again being the average size of $f(m)$.

The index of a matrix $R=R_q$ with $q\in\KK$ is given by
\begin{equation}
   \varSigma (R_q ) \; = \;
   \N \bigl( \lvert q / \mathrm{cont}^{}_{\KK} (q) \rvert^2 \bigr).
\end{equation}
It is thus of the form $\N \bigl( k +\ell\sqrt{2}\ts \bigr)= k^2 - 2
\ell^2$ with $k,\ell\in\ZZ$.  Since each prime $\pi\in\ZZ[\sqrt{2}\,]$
has a decomposition $\pi=p\bar{p}$ in $\KK$, one sees that all
indices of this form indeed occur.  This gives
\begin{coro}
  The elementary coincidence spectrum of the\/ $\ZZ[\sqrt{2}\,]$-module\/
  $\I (\KK)$ is
\[
   \begin{split}
   \varSigma \bigl(\mathrm{SO} (3,\QQ(\sqrt{2}\,))\bigr)
   & \; = \; \{ m\in\NN \mid m= k^2 - 2 \ell^2 \, ,
   \mbox{ with } k,\ell\in\ZZ \} \\[2mm]
   & \; = \;
   \{ 1,2,4,7,8,9,14,16,17,18,23,25, \ldots \}\, .
   \end{split}
\]
   The spectrum is a monoid and consists of all natural numbers that
   have all primes congruent to\/ $\pm 3$ $(\mbox{mod } 8)$ in their
   factorizations occuring with an even power. \qed
\end{coro}

A related module in $3$-space is $\ZZ[\sqrt{2}\,]^3$, realized as $\I
(\LL)$ with the order $\LL$ of \eqref{lattice}. The treatment of this
module parallels the situation already met with the icosian ring, see
above.

\section{Extensions and outlook}

The picture for the coincidence site modules in 3-space that emerge
from single intersections seems to be rather clear now. Still, there
are several important questions that have not been addressed here at
all. First of all, having counted the CSMs, one would like to have
a finer classification into Bravais types. Some results in this
direction, for the cubic lattices, are presented in \cite{Z1}.
Then, there is no compelling reason to stop at single intersections,
and recent developments make an extension to multiple coincidence
site modules desirable. While this is rather well understood in the planar
case \cite{BG}, by algebraic means based upon the arithmetic of
cyclotomic fields, only first steps exist in $3$-space, see
\cite{Z2,Z3} for some results.

As we briefly mentioned in passing, several of our above findings can
be generalized beyond the class number $1$ situation, and we hope to
report on that soon \cite{BPnew}. Quaternions are certainly also
helpful in $4$-space, and further progress in this direction is in
sight \cite{BZ}, at least for simple coincidences.  More complicated,
however, seems the situation in higher dimensions, even for the class
of root lattices, and that might be a good problem to solve.

\bigskip
\section*{Appendix}

In the proof of Lemma~\ref{ideal-map}, we needed a global 
version of a uniqueness result on intersecting maximal orders in
quaternion division algebras, which is to some extent implicit in
\cite{R,V}, but which we could not find in explicit form and which
does not seem to have an analogue in the more general context of
central simple algebras \cite{R,Al}.  For completeness, we sketch
here a localization argument for this \cite{St2,Al} that was
communicated to us by U.\ Staemmler. For comparison with other
work, we recall that the intersection of two (not necessarily 
distinct) maximal orders is called an \emph{Eichler order}.

\begin{lemma}\label{appendix} \cite{St2}
  Let\/ $\HH (K)$ be the quaternion division algebra over the
  real algebraic number field\/ $K$, and let\/ $M$, $N$ and $O$
  be three maximal orders of\/ $\HH (K)$ that satisfy\/
  $O\cap M = O\cap N$. Then $M=N$.
\end{lemma}

\begin{proof}
We write $H=\HH(K)$ for this proof.  One has $M=N$ if and only if
all local completions of $M$ and $N$ at prime ideals $\mathfrak p$
of $K$ are equal, i.e., $M_{\mathfrak p}=N_{\mathfrak p}$, see
\cite[Prop.~III.5.1]{V} and \cite[Sec.~3]{R} for details in
the present setting.  Note that $E:= O\cap M = O\cap N$ is an
Eichler order.  Local completion preserves intersections and the
property of being a maximal order, see \cite[p.~84]{V}.  In
particular, $E_{\mathfrak p} = O_{\mathfrak p} \cap M_{\mathfrak p}
= O_{\mathfrak p} \cap N_{\mathfrak p}$ remains an Eichler order in
$H_{\mathfrak p}$.

Either $H_{\mathfrak p}$ is a skew field (when $\mathfrak{p}$ divides
the discriminant of $H$ over $K$) or else it is isomorphic to the
matrix ring $\mathrm{Mat} (2,K_{\mathfrak p})$ by Wedderburn's
theorem, see \cite[Sec.~52]{O}. In the latter case, $M_{\mathfrak p} =
N_{\mathfrak p}$ by \cite[Lemme~II.2.4]{V}, which states that Eichler
orders uniquely specify the two intersecting maximal orders here.  In
the former case, $M_{\mathfrak p} = N_{\mathfrak p}$ because
$H_{\mathfrak p}$ has a unique maximal order by
\cite[Lemme~II.1.5]{V}.  So, $M_{\mathfrak p} = N_{\mathfrak p}$ for
all prime ideals $\mathfrak p$ and hence $M=N$.
\end{proof}

Note that this proof actually shows a slightly stronger statement:
if $L,M,N,O$ are maximal orders with $L\cap M=N\cap O$, one has
$\{L,M\}=\{N,O\}$. But the lemma as stated is enough for our purposes.

\bigskip
\section*{Acknowledgements}

It is a pleasure to express our gratitude to Ute Staemmler, Alfred
Weiss and Peter Zeiner for their cooperation and for useful hints on
the manuscript.  Interesting discussions with Herbert Abels, Robert
V.\ Moody and Claus-Michael Ringel are gratefully acknowledged.
M.~B.\ would like to thank the University of Queensland for support in
form of a Raybould Visiting Fellowship and the Department of
Mathematics for hospitality, where a substantial part of this work was
done. It was also supported by the German Research Council (DFG),
within the CRC 701.

\end{document}